\documentclass[aos]{imsart}

\RequirePackage{amsthm,amsmath,amsfonts,amssymb}
\RequirePackage[numbers,sort&compress]{natbib}
\RequirePackage[colorlinks,citecolor=blue,urlcolor=blue]{hyperref}

\startlocaldefs
\theoremstyle{plain}

\newtheorem{Theorem}{Theorem}[section]

\newtheorem{Proposition}[Theorem]{Proposition}
\newtheorem{Corollary}[Theorem]{Corollary}
\theoremstyle{remark}
\newtheorem{Assumption}[Theorem]{Assumption}
\newtheorem{Definition}[Theorem]{Definition}
\newtheorem{Remark}{Remark}

\usepackage[utf8]{inputenc}
\usepackage{amsfonts}
\usepackage{dsfont}
\usepackage{amsgen,amsmath,amstext,amsbsy,amsopn,amssymb,subcaption}
\usepackage{graphicx}

\usepackage[shortlabels]{enumitem}
\usepackage{natbib}
\usepackage{booktabs}
\usepackage{amsthm}
\usepackage{hyperref}
\usepackage{blkarray}
\usepackage{algorithm}
\usepackage{algorithmic}
\usepackage{url}
\usepackage{longtable}
\usepackage{comment}
\usepackage{stmaryrd}

\SetSymbolFont{stmry}{bold}{U}{stmry}{m}{n}
\DeclareMathAlphabet\mathbfcal{OMS}{cmsy}{b}{n}

\newcommand{\dist}{{\rm dist}}

\makeatletter
\newcommand*{\rom}[1]{\expandafter\@slowromancap\romannumeral #1@}
\makeatother

\usepackage[dvipsnames]{xcolor}

\usepackage{stmaryrd}
\SetSymbolFont{stmry}{bold}{U}{stmry}{m}{n}
\usepackage{bm}

\DeclareMathAlphabet\mathbfcal{OMS}{cmsy}{b}{n}

\newcommand{\abs}[1]{\left| #1 \right|}
\newcommand{\norm}[1]{\left\| #1 \right\|}

\usepackage[dvipsnames]{xcolor}

\newcommand{\argmax}{\mathop{\rm arg\max}}

\newcommand{\spa}{{\rm span} }
\newcommand{\Proj}{{\rm P} }
\newcommand{\poly}{\operatorname{poly}}



\newcommand{\bI}{\mathbf{I}}

\newcommand{\cB}{\mathcal{B}}

\newcommand{\cD}{\mathcal{D}}

\newcommand{\cH}{\mathcal{H}}

\newcommand{\cK}{\mathcal{K}}

\newcommand{\cM}{\mathcal{M}}
\newcommand{\cN}{\mathcal{N}}

\newcommand{\cP}{\mathcal{P}}
\newcommand{\cQ}{\mathcal{Q}}

\newcommand{\cS}{\mathcal{S}}
\newcommand{\cT}{\mathcal{T}}

\newcommand{\cV}{\mathcal{V}}

\newcommand{\cX}{\mathcal{X}}

\def\bbB{\mathbb{B}}
\def\R{\mathbb{R}}

\def\bbZ{\mathbb{Z}}

\newcommand{\gap}{{\mathsf{gap} }}

\newcommand{\sfm}{{\mathsf{m}}}

\newcommand{\sfH}{{\mathsf{H}}}

\newcommand{\opt}{{\mathsf{opt}}}
\newcommand{\PP}{\mathbb{P}}
\newcommand{\E}{\mathbb{E}}

\newcommand{\Reg}{\mathsf{Reg} }
\newcommand{\Leb}{\mathsf{Leb} }
\newcommand{\sgn}{\mathsf{sgn} }
\newcommand{\est}{\mathsf{est}}
\newcommand{\idc}{\mathbb{I}}

\newcommand{\lp}{\underline{p} }

\newcommand{\lmu}{\underline{\mu}}
\newcommand{\umu}{\overline{\mu}}
\newcommand{\betamin}{\underline{\beta} }
\newcommand{\betamax}{\overline{\beta} }
\newcommand{\Cov}{{\rm{Cov}}}
\newcommand{\sff}{{\mathsf{f}}}
\newcommand{\cf}{c_{\sff}}
\newcommand{\sfX}{{\mathsf{X}}}
\newcommand{\sfc}{{\mathsf{c}}}

\endlocaldefs

\begin{document}

\begin{frontmatter}
\title{Nonparametric Bandits with Single-Index Rewards:
	\\
	Optimality and Adaptivity}
\runtitle{Nonparametric Bandits with Single-Index Rewards}

\begin{aug}
\author[A]{\fnms{Wanteng}~\snm{Ma} \ead[label=e1]{wanteng@wharton.upenn.edu}}
\and
\author[B]{\fnms{T. Tony}~\snm{Cai}\ead[label=e2]{tcai@wharton.upenn.edu}}

\address[A]{Department of Statistics and Data Science, University of Pennsylvania \printead[presep={,\ }]{e1}}

\address[B]{Department of Statistics and Data Science, University of Pennsylvania \printead[presep={,\ }]{e2}}
\end{aug}

\begin{abstract}
	Contextual bandits are a central framework for sequential decision-making, with applications ranging from recommendation systems to clinical trials. While nonparametric methods can flexibly model complex reward structures, they suffer  from the curse of dimensionality. We address this challenge using a single-index model, which projects high-dimensional covariates onto a one-dimensional subspace while preserving nonparametric flexibility.

	We first develop a nonasymptotic theory for offline single-index regression under a rank-monotone response, combining maximum rank correlation for index estimation with local polynomial regression. Building on it, we propose an epoch-based batched bandit algorithm and prove upper and lower bounds on the regret over a structured class of single-index instances, which match up to logarithmic factors over a broad range of smoothness and margin parameters. In this range the minimax regret is of order $ n^{1-\frac{ (1+\alpha)\beta }{2\beta+1 } } + n\left(\frac{d}{n}\right)^{\frac{1+\alpha}{2}}$, where $\beta$ is the smoothness of the link and $\alpha$ is the margin parameter. Thus the ambient dimension $d$ enters only through a parametric term, with a phase transition at $d\asymp n^{\frac{1}{2\beta+1}}$, and in fixed dimension the regret matches the optimal rate of one-dimensional nonparametric bandits. We further show that adaptation to an unknown smoothness level is impossible in general, but is achievable under a standard self-similarity condition.

\end{abstract}



\end{frontmatter}

\section{Introduction}
\label{sec:intro}

Online decision-making has become central to modern data-driven systems, underpinning applications ranging from recommendation and advertising to clinical trials and resource allocation. A key objective in these settings is reward maximization under uncertainty, where an agent learns optimal decisions through feedback over time. The \textit{bandit problem} provides a foundational framework for this challenge and has been extensively studied across multiple variants, including multi-armed bandits \citep{lai1985asymptotically,perchet2013multi}, linear contextual bandits \citep{chu2011contextual,abbasi2011improved}, and adversarial bandits \citep{auer2002nonstochastic,bubeck2012regret}.

Among these frameworks, \textit{contextual bandits} stand out for their ability to incorporate side information, enabling personalized and adaptive decision-making. They have achieved remarkable success in recommendation systems \citep{li2010contextual}, adaptive clinical trials \citep{durand2018contextual}, online advertising \citep{chapelle2011empirical}, and intelligent tutoring systems \citep{rafferty2011faster}. While many studies assume linear reward structures for tractability \citep{goldenshluger2013linear,bastani2020online}, real-world feedback is often nonlinear: in the Netflix artwork selection problem \citep{chandrashekar2017artwork,gur2022smoothness}, for instance, a user's willingness to play a movie title is a non-linear function of the user's preference over two artworks portraying the title. This motivates \textit{nonparametric contextual bandits}, which allow greater modeling flexibility but suffer from the \textit{curse of dimensionality}: for the algorithms studied in the literature \citep{rigollet2010nonparametric,perchet2013multi,gur2022smoothness}, the sample complexity required to learn a good policy grows exponentially with the dimension $d$, making these methods challenging to deploy in practice.

Motivated by these challenges, we consider \textit{single-index bandits}, where each reward function depends on covariates only through an unknown one-dimensional projection, called an \textit{index}, while retaining nonparametric flexibility in the link function. By reducing the effective dimension to one, single-index contextual bandits preserve the expressive power of nonparametric models while alleviating the curse of dimensionality, enabling efficient learning in complex, high-dimensional decision-making problems. The \textit{single-index model} has been extensively studied in nonparametric regression and widely applied in practice \citep{ichimura1993semiparametric,horowitz1998semiparametric,xia2008multiple}. We show that over a structured class of single-index bandit instances, specified in Sections~\ref{sec:set-up}\textendash\ref{sec:bandits} by responses that are monotone in the index, a nondegenerate index-projected design and a learnable margin, it is possible to design an efficient batched algorithm whose regret is minimax optimal up to logarithmic factors in the range of smoothness and margin parameters covered by our lower bounds, and that there the ambient dimension $d$ affects the optimal regret only through an additive parametric term. In particular, in the fixed-dimensional regime the optimal regret coincides with the one-dimensional nonparametric bandit rate, thereby achieving effective dimension reduction and overcoming a key barrier in high-dimensional online learning.

We briefly describe the formulation of the multi-armed contextual bandit problem as follows. At each time step $t = 1, \ldots, n$, a decision maker is presented with a covariate vector (or context) $X_t \in \mathbb{R}^d$ and selects one arm from a $K$-armed bandit $\{Y_t^{(k)}\}_{k=1}^K$ ($K \ge 2$). The decision maker then observes a reward $Y_t = Y_t^{(\pi_t)}$ corresponding to the chosen arm $\pi_t \in [K]$. The term ``contextual'' indicates that the decision at time~$t$ is made based on the contextual information~$X_t$, which is observed prior to choosing the arm.
Suppose further that the $d$-dimensional contexts $X_t$, together with the rewards for the $K$ arms $\{Y_t^{(k)}\}_{k=1}^K$, are independent and identically distributed across time steps $t = 1, 2, \ldots, n$.
\begin{equation}\label{eq:bandit-dist}
	\left(X_t, Y_t^{(1)}, Y_t^{(2)},\dots,Y_t^{(K)}\right)\sim P_{X, \{Y^{(k)}\}_{k=1}^K}.
\end{equation}
Define the (expected) reward function for each arm as
\begin{equation*}
	g_k(X_t):= \E [Y^{(k)}_t|X_t].
\end{equation*}
Our goal is to construct a policy $\pi=\{\pi_t\}_{t=1}^n$ that minimize the cumulative expected regret with respect to the optimal policy $\pi^\star(X_t)=\argmax_{k\in{K}} \{g_k(X_t)\}$, i.e.,
\begin{equation}\label{eq:regret}
	\Reg(\pi) = \E\left[\sum_{t=1}^n \left(Y^{(\pi^\star(X_t) )}_t-  Y^{(\pi_t)}_t \right)\right] = \E\left[\sum_{t=1}^n  \left( g_{\pi^\star(X_t)}(X_t) - g_{\pi_t}(X_t)\right)\right].
\end{equation}
In this paper, we focus on the setting where the reward functions $g_k(X_t)$ follow a single-index model (with the detailed definition provided in Section~\ref{sec:set-up}) for dimension reduction. We study the minimax regret for single-index bandits and develop online algorithms whose regret attains it, up to logarithmic factors, in the range covered by our lower bounds.

\subsection{Main contribution}

Our contribution is three-fold:
\begin{itemize}
	\item \textit{Minimax regret for single-index bandits.} We first propose an offline single-index regression method based on an index plug-in approach and establish a sharp non-asymptotic convergence rate. For responses that are monotone in the index and links with H\"older smoothness $\beta$, we employ maximum-rank correlation \citep{han1987non,fan2020rank} for index estimation and show that it attains $\ell_2$ accuracy of order $\widetilde{O}(\sqrt{d/n})$ with overwhelming probability. On the way we identify a set of primitive conditions (a nondegenerate monotone response, a conditional density floor on the index, and a nondegenerate design covariance) under which the population rank-correlation objective is locally strongly concave with explicit constants. The estimated index is then plugged into a one-dimensional local polynomial regression. Building on this, we design an epoch-based batched bandit algorithm with $O(\log n)$ data-driven epochs and a two-step arm elimination: arms are first preselected using the previous epoch's estimates, then refined using the most recent data. We prove a regret bound of order $\widetilde{O}(n^{1-\frac{(\alpha+1)\beta}{2\beta+1}})$, where $\alpha$ is the margin parameter, and a minimax lower bound of the same order over the structured class above for a broad range of $(\alpha,\beta)$; in fixed dimension this rate is independent of $d$, demonstrating effective dimension reduction.

	\item \textit{Adaptivity to unknown smoothness.} We next investigate adaptation when the smoothness level $\beta$ is unknown. We show that, analogous to the classical nonparametric bandit setting, adaptation to an unknown $\beta$ is impossible without additional structural assumptions \citep{gur2022smoothness}. However, under a standard self-similarity condition, adaptation becomes feasible by first estimating $\beta$ prior to regret minimization. To this end, we propose a link-smoothness estimation procedure with index plug-in, which yields an ``under-smoothed'' estimator of $\beta$ with accuracy $O\left(\frac{\log\log n}{\log n}\right)$. This estimation step enables our algorithm to attain the optimal regret up to logarithmic factors even when $\beta$ is unknown, for the ranges of $\beta$ and $d$ specified in Theorem~\ref{thm:adaptive-regret-ub}.

	\item \textit{Single-index bandits in increasing dimensions.} Finally, we study the minimax regret when the dimension $d$ grows with the sample size $n$ (with $d\ll n$). We show that, in the increasing-dimension regime, the regret of our algorithm is of order $n\left(\left(\frac{ 1 }{n}\right)^{\frac{ (1+\alpha)\beta }{2\beta+1 } } +   \left(\frac{d}{n}\right)^{\frac{1+\alpha}{2}}\right)$ up to logarithmic factors, and that this is the minimax rate over a range of $(\alpha,\beta)$ given in Section~\ref{sec:increasing-dim}. It reveals a two-phase phenomenon: a nonparametric phase $n^{1-\frac{(\alpha+1)\beta}{2\beta+1}}$, corresponding to the optimal regret of one-dimensional nonparametric bandits, and a parametric phase $n\left(\frac{d}{n}\right)^{\frac{1+\alpha}{2}}$, corresponding to the optimal regret of linear bandits \citep{bastani2021mostly,mao2025optimal}. Our algorithm attains this rate without knowing which phase it is in, thereby bridging the two.
\end{itemize}

It is worth noting at the outset that our optimality is stated over a mild structured subclass of single-index bandits rather than over all of them. Both the upper and the lower bounds are established over the class $\cP$ of Definition~\ref{def:dist-bandit} with link functions in $\cH(\beta,L)$, and they match only in the ranges of $(\alpha,\beta)$ stated in Theorems~\ref{thm:minimax-lb} and~\ref{thm:minimax-lb-phase}; Section~\ref{sec:discuss} records the range that remains open. The class requires responses that are monotone in the index in the rank sense of Assumption~\ref{asp:model}, regularity of the index-projected design, nondegeneracy of the rank-correlation criterion, and learnable optimal and margin regions for each arm. The first three make the index estimable at a parametric rate under bandit sampling, and the last makes the link estimable where decisions are difficult; each is discussed where it is used, in Sections~\ref{sec:set-up}\textendash\ref{sec:bandits}. Our policy is also \emph{batched}, updating only at the end of each of $M=O(\log n)$ epochs.

\subsection{Related literature}

Bandit algorithms have been studied extensively; see \cite{lattimore2020bandit} for an overview. Closest to our setting are linear and generalized linear bandits, where the reward is a known function of a linear score \citep{chu2011contextual,li2017provably}, and \emph{nonparametric bandits}, pioneered by \cite{rigollet2010nonparametric} and \cite{perchet2013multi}, in which no parametric form is imposed. The latter partition the covariate space into bins: \cite{perchet2013multi} proposed ABSE, a successive-arm-elimination method that is minimax optimal for $\beta\le 1$, and \cite{hu2022smooth} obtained optimal smooth regret rates for $\beta\ge 1$ by exploiting information within bins. All of these bin-based bounds scale exponentially with $d$, typically as $O((\poly(\beta))^{d})$, which is what makes dimension reduction valuable. Two routes to it have been pursued: sparsity, where the score depends on $s\ll d$ coordinates and best subset selection yields regret nearly free of the ambient dimension \citep{chen2024nearly}; and the single index, which we follow, keeping all $d$ coordinates but replacing the parametric score by an unknown link. The single-index model is the classical device for it: the minimax MISE rate for nonparametric regression, $O(n^{-\frac{2\beta}{2\beta+d}})$, improves to $O(n^{-\frac{2\beta}{2\beta+1}})$ under a single index \citep{gaiffas2007optimal}. Section~C of the Supplementary Material~\citep{suppSIB} reviews both literatures in more detail.

\paragraph*{Semiparametric and single-index bandits} Two recent works are closely related to ours. \cite{arya2025semi} study batched multi-armed bandits with covariates under a semiparametric single-index reward that is \emph{shared} across arms, $g_k(x)=f(v^\top x)+ \theta_k$, estimating the nonparametric component by binning; their focus is the \textit{global} (fixed-arm) objective with a small number of batches, and their guarantees are stated for Lipschitz-type links. \cite{kang2025single} study generalized linear contextual bandits with an unknown Lipschitz link, again with a single index shared by all arms, and additionally require the covariate distribution to be known. Both show that rank-based or isotonic estimation of an unknown monotone link is compatible with sequential decision-making, and neither is subsumed by our results. Relative to them, our contribution differs in four respects. First, we allow \emph{arm-specific} indices $v_1,\ldots,v_K$, so that the optimal decision regions are determined by $K$ distinct projections and cannot be reduced to a single scalar summary of the context; this is the source of Assumption~\ref{asp:decision-region} and of the two-step elimination in Algorithm~\ref{alg:sim-bandit}. Second, we treat general H\"older smoothness $\beta>0$ and establish matching minimax lower bounds over a broad range of smoothness and margin parameters. Third, our algorithm uses no bin partition, and hence avoids the $(\log^{1/\beta}n \vee 2)^d$ factor carried by binning-based methods. Fourth, we address adaptation to unknown smoothness and the increasing-dimension phase transition, and we do not assume knowledge of $P_X$. Table~\ref{tab:comparison} summarizes the comparison.

Single-index structure has also been used in dynamic pricing, where the action is a continuous price that enters the index \citep{fan2024policy,luo2024distribution,wang2025tight,han2026semiparametric}; the conditions there are specific to pricing and neither imply nor are implied by ours (Section~C of the Supplementary Material~\citep{suppSIB}).

\begin{table}[h!]
	\centering
	\small
	\begin{tabular}{lccccc}
		\toprule
		& reward model & link & index & bin partition & lower bound \\
		\midrule
		\cite{perchet2013multi} & nonparametric & $\beta\le 1$ & \textemdash & yes & yes \\
		\cite{hu2022smooth} & nonparametric & $\beta>0$ & \textemdash & yes & yes \\
		\cite{li2017provably} & generalized linear & known & arm-specific & no & no \\
		\cite{arya2025semi} & single-index & monotone, Lipschitz & shared & yes & no \\
		\cite{kang2025single} & single-index & monotone, Lipschitz & shared & no & no \\
		This paper & single-index & monotone, $\beta>0$ & arm-specific & no & yes \\
		\bottomrule
	\end{tabular}
	\caption{Comparison with the relevant literature. ``Index'' indicates whether the low-dimensional projection is common to all arms or arm-specific. \cite{kang2025single} additionally assume that the covariate distribution $P_X$ is known.}
	\label{tab:comparison}
\end{table}

Adaptation is another important topic in nonparametric bandits. Adaptation to unknown smoothness is generally impossible without further assumptions \citep{locatelli2018adaptivity,cai2022stochastic,gur2022smoothness}, but becomes possible under a self-similarity condition \citep{picard2000adaptive,gine2010confidence}, an idea extended to nonparametric bandits via Lepski's method \citep{lepski1997optimal}. Whether such adaptation remains feasible under an additional single-index structure has remained open.

\subsection{Organization}

The rest of the paper is organized as follows. Section~\ref{sec:set-up} introduces a formal mathematical formulation of the nonparametric bandit problem with single-index rewards. Section~\ref{sec:regression} presents our (offline) single-index regression method and establishes sharp non-asymptotic error bounds. Building on these results, Section~\ref{sec:bandits} develops our online algorithm for single-index bandits. Section~\ref{sec:opt-rates} bounds the regret of our algorithm from above and the minimax regret from below, and shows that the two bounds match in the range stated there. Section~\ref{sec:adaptation} investigates adaptation to unknown smoothness levels. Section~\ref{sec:increasing-dim} extends the problem to increasing dimensions and examines the resulting phase transition in minimax optimal rates. Finally, Section~\ref{sec:simulations} presents numerical experiments that corroborate our theoretical findings and compares the proposed algorithm with nonparametric, parametric and neural alternatives on simulated and on benchmark data. The proofs and the remaining experiments are given in Sections~A and~F of the Supplementary Material~\citep{suppSIB}.

\subsection{Notation}

The following notation will be used throughout the paper. For a function $f$, we use $f^{(i)}$ to denote its $i$-th derivative. For a random variable $Y$, $Y^{(k)}$ represents the variable corresponding to the $k$-th arm. We also use $g_{(i)}(x)$ to denote the $i$-th largest value among the functions $\{g_k\}_{k\in[K]}$ evaluated at $x$, with the precise definition provided in the next section. For parameters associated with each arm that share a common notation, such as $v_k$, we use $v_{k,i}$ to represent the $i$-th entry of the parameter in arm $k$. For vectors without an arm specification, such as $x$, $x_i$ denotes the $i$-th entry of $x$.
For $\beta>0$ we write $\lfloor \beta \rfloor$ for the largest integer strictly smaller than $\beta$, which is also the local polynomial degree used throughout. The notation $\norm{\cdot}_2$ represents the $\ell_2$ norm. For two scalars $a,b$, we define $a \vee b = \max\{a,b\}$ and $a \wedge b = \min\{a,b\}$. The $\ell_2$ ball is denoted by $\bbB_2(x,r)=\left\{x'\in\R^d \mid \norm{x'-x}_2\le r \right\}$. For any positive integer $K$, we let $[K] = \{1,2,\ldots,K\}$ denote the corresponding index set. For a set $\cX \subseteq \R^d$, its projection onto a direction $v$ is defined as $v \circ \cX = {v^\top x \mid x \in \cX}$. For a scalar $x$ and $r > 0$, the $r$-neighborhood of $x$ is denoted by $U(x, r) = [x - r, x + r]$. For an interval or general set $I \subseteq \R$, the $r$-neighborhood of $I$ is defined as $U(I, r) =\{x \mid \dist(x, I) \le r\}$, where $\dist(x, I)$ represents the distance from $x$ to the closest element in $I$. We use the standard asymptotic notation $O(\cdot)$, with the same meaning as $\lesssim$, and write $f(n)=\widetilde{O}(g(n))$ when $f(n)=O(g(n)\poly(\log n))$.

\section{Nonparametric Bandits with Single-Index Rewards}
\label{sec:set-up}

We introduce in this section a formal framework for the nonparametric bandit problem with single-index rewards.
Given the unknown bandit distribution in \eqref{eq:bandit-dist}, our goal is to minimize the cumulative regret defined in \eqref{eq:regret}.
We denote the marginal distribution of $X_t$ by $P_X$, with support $\cX$, and assume that it satisfies certain regularity conditions (to be specified later).
As is standard in bandit problems, once arm $k \in [K]$ is selected at time $t$, only $Y_t^{(k)}$ is observed (i.e., $Y_t = Y_t^{(k)}$), while the rewards for all other arms at time $t$ remain unobserved.
Under this partial-feedback structure, we aim to design a policy $\pi = \{\pi_t\}_{t=1}^n$ for arm selection, where each $\pi_t$ is an \textit{admissible} randomized decision rule $\pi_t:\R^d \to [K]$ that chooses which arm to pull at time $t$ based on $X_t$, and depends only on observations strictly prior to time $t$.
Throughout our analysis, we assume that the number of arms remains at a constant level, i.e., $K = O(1)$.

Without imposing additional structural assumptions, minimizing \eqref{eq:regret} is challenging due to the large functional space.
In this paper, we assume that for each arm $k \in [K]$, the expected reward follows a single-index model \citep{ichimura1991semiparametric,manski1988identification}, i.e.,
\begin{equation}\label{eq:sim}
	g_k(X_t) = \E [Y^{(k)}_t \mid X_t] = f_k(v_k^\top X_t),
\end{equation}
where the \textit{link} function $f_k:\R\to[0,1]$ and the \textit{index} vector $v_k \in \cV \subseteq \R^d$ are both unknown to the decision-maker. The range $[0,1]$ is fixed only for definiteness, so that the Bernoulli response model below is covered; all results hold for links with values in any fixed bounded interval after rescaling.
We assume that $\cX$ and the parameter space $\cV$ are compact; the conditions imposed on the design are recorded in \eqref{eq:design-bounded} below.

However, due to the flexible choices of $f_k$ and $v_k$, the model \eqref{eq:sim} is not identifiable without additional constraints.
To ensure identifiability of both the link functions and index vectors, we impose the normalization $v_{k,1} = 1$ for each arm $k$.
For further discussion on identifiability in single-index models, see, e.g., \cite{ichimura1993semiparametric,manski1988identification}, and Chapter~2 of \cite{horowitz2009semiparametric}.

For each arm $k\in[K]$, we pose the reward in a regression form:
\begin{equation}\label{eq:reg-form}
	Y^{(k)}_t = g_k(X_t)+ \xi_{t}^{(k)} =  f_k(v_k^\top X_t) +\xi_{t}^{(k)},
\end{equation}
where $\xi_{t}^{(k)}:=Y^{(k)}_t-f_k(v_k^\top X_t)$ is the regression residual, which may be $X_t$-dependent. We emphasize that $\xi_{t}^{(k)}$ is an \textit{induced} residual rather than a primitive error term: the primitive object is the conditional law of $Y^{(k)}_t$ given $X_t$, which we now restrict. The following assumption collects all our modeling conditions on the response.

\begin{Assumption}[Single-index rank-monotone response]\label{asp:model}
	For each arm $k\in[K]$ and each $t$, the conditional law of $Y_t^{(k)}$ given $X_t$ depends on $X_t$ only through the index $v_k^\top X_t$, and it is \emph{rank-monotone} in the index: if $(X,Y)$ and $(X',Y')$ are independent copies of $(X_t,Y^{(k)}_t)$, then
	\begin{equation}\label{eq:stoch-monotone}
		\PP\left(Y > Y' \mid v_k^\top X = u,\ v_k^\top X'=u'\right)\ \ge\ \PP\left(Y < Y' \mid v_k^\top X = u,\ v_k^\top X'=u'\right) \quad \forall u\ge u' .
	\end{equation}
	Its conditional mean is $\E[Y^{(k)}_t\mid X_t]= f_k(v_k^\top X_t)$ with $f_k\in\cH(\beta,L)$ as in Definition~\ref{def:smooth-class}, so that $\E[\xi_{t}^{(k)}\mid X_t]=0$, and $\xi_{t}^{(k)}$ is conditionally sub-Gaussian: $\E\left[\exp\left((\xi_{t}^{(k)})^2/C\right) \mid X_t\right]\le 2$ for some constant $C=O(1)$.
\end{Assumption}

Assumption~\ref{asp:model} isolates the two distinct roles played by the response model in our analysis. The rank monotonicity~\eqref{eq:stoch-monotone} is used \emph{only} for index estimation: it is what guarantees that the population rank correlation is maximized at $v=v_k$, and it is what justifies the use of the MRC estimator in Section~\ref{sec:regression}. It is implied by stochastic monotonicity of the conditional law in the index, as in the generalized regression model of \cite{han1987non}, but it constrains the responses only through their ranks and does not force the mean link $f_k$ to be monotone (Section~G.5 of the Supplementary Material~\citep{suppSIB}); the link class below therefore carries no monotonicity constraint. Notice that the residual $\xi^{(k)}_t$ itself does not enter the analysis of the MRC step, since the criterion~\eqref{eq:mrc} depends on the responses only through their ranks; the response law enters that step only through the concordance function of Assumption~\ref{asp:nondeg:1}, and $\xi^{(k)}_t$ enters the analysis only through the conditional mean-zero and sub-Gaussian conditions in the local polynomial step.

A convenient sufficient condition for \eqref{eq:stoch-monotone}, under which $f_k$ is also nondecreasing, is the monotone transformation representation
\begin{equation}\label{eq:D-rep}
	Y^{(k)}_t = D\left(f_k(v_k^\top X_t)+u_t^{(k)}\right),
\end{equation}
where $u_t^{(k)}$ is independent of $X_t$ and $D$ is nondecreasing and \textit{calibrated}, in the sense that $\E D(z+u^{(k)}_t) = z$ for all $z$ in the range of $f_k$; calibration is what delivers $\E[Y_t^{(k)}\mid X_t]=f_k(v_k^\top X_t)$. Typical examples include general i.i.d.\ noise regression, where $\xi_{t}^{(k)}$ is independent of $X_t$ and $D(z)=z$; and Bernoulli response regression $Y^{(k)}_t\mid X_t \sim\operatorname{Ber}( f_k(v_k^\top X_t))$, which corresponds to $D(z)=\idc\{z>\frac{1}{2}\}$ with $u_t^{(k)}\sim \operatorname{Unif}(-\frac{1}{2},\frac{1}{2})$. More generally, \eqref{eq:D-rep} covers monotone transformation models in the sense of \cite{han1987non}, while Assumption~\ref{asp:model} also covers response models that need not be of the form \eqref{eq:D-rep}, for instance averaged binomial responses $Y_t^{(k)}=m^{-1}\operatorname{Bin}\left(m, f_k(v_k^\top X_t)\right)$ with $m\ge 2$ fixed, whose conditional law is stochastically increasing in the index and sub-Gaussian, but is not of the form \eqref{eq:D-rep} with a single noise variable $u^{(k)}_t$ independent of $X_t$.

To characterize $f_k$, we introduce the standard Hölder smoothness class as follows:

\begin{Definition}[Smoothness class]\label{def:smooth-class} We say that the link function $f$ belongs to the Hölder class $\cH_0(\beta, L)$ for $\beta>0$ and $L>0$, which is the collection of functions $f$: $\R\rightarrow [0,1]$ that are $\lfloor \beta \rfloor$ times continuously differentiable and for any $x, x^{\prime} \in \R$ satisfy the following, where $\lfloor \beta \rfloor$ denotes the largest integer strictly smaller than $\beta$:
	$$
	\left|f^{(\lfloor \beta \rfloor )}\left(x\right)-f^{(\lfloor \beta \rfloor)}\left(x^{\prime}\right)\right| \leq L \abs{x-x^{\prime}}^{\beta -\lfloor \beta \rfloor }.
	$$
	Define our function class of interest as
	\begin{equation}\label{eq:monotone}
		\mathcal{H}(\beta, L)=
		\begin{cases}\mathcal{H}_0(\beta, L), & \beta<1, \\
			\mathcal{H}_0(\beta, L) \cap \mathcal{H}_0(1, L), & \beta \geq 1 .
		\end{cases}
	\end{equation}
\end{Definition}
Here, $\mathcal{H}(\beta, L)$ is defined based on $\cH_0(\beta, L)$ for minimax purposes. The Hölder smooth class for multivariate functions can be similarly defined using Taylor polynomials. See, e.g, \cite{audibert2007fast}. For $\beta\ge1$, $\cH(\beta,L)$ also imposes global Lipschitz continuity with constant $L$, as is standard for smooth nonparametric bandits \citep{hu2022smooth,gur2022smoothness}; it converts an index estimation error into a reward estimation error, $\abs{f(\widehat{v}^\top x)-f(v^\top x)}\le L\norm{\widehat{v}-v}_2\norm{x}_2$.

Finally, we record the conditions imposed on the covariate distribution. Throughout Sections~\ref{sec:set-up}\textendash\ref{sec:adaptation}, where the dimension $d$ is fixed, we assume that the design is bounded and that the index range is bounded for simplicity:
\begin{equation}\label{eq:design-bounded}
	\sup_{x\in\cX}\norm{x}_2\le C_{\cX}, \qquad \sup_{x\in\cX}\abs{v_k^\top x}\le C_I \quad \text{ for every } k\in[K],
\end{equation}
for constants $C_{\cX},C_I>0$. We use Condition~\eqref{eq:design-bounded} through the resulting sub-Gaussianity of the one-dimensional projections, $\sup_{u\in\bbB_2(1)}\norm{u^\top X}_{\psi_2}\le 2C_{\cX}$. In Section~\ref{sec:increasing-dim}, where $d$ grows with $n$, boundedness of $\norm{X}_2$ by a $d$-free constant is incompatible (Section~D.3 of the Supplementary Material~\citep{suppSIB}), and the conditions needed there are stated separately as Assumption~\ref{asp:normalization}.
To ease the notation, in what follows, we shall write the optimal arm and sub-optimal arm as
\begin{equation}\label{eq:margin-def}
	\begin{aligned}
		g_{(1)}(x)& =\max_{k\in[K]} \left\{g_{k}(x) \right\},\\
		g_{(2)}(x) & =\max_{k\in[K]} \left\{g_{k}(x)\mid g_{k }(x)<g_{(1)} (x)\right\},
	\end{aligned}
\end{equation}
or   $g_{(2)}(x)  =g_{(1)}(x)$ if   $\left\{k\mid g_{k}(x)<g_{(1)} (x)\right\} = \emptyset$.

\section{A Non-Asymptotic Theory for Single-Index Regression}
\label{sec:regression}

In this section, we study the offline regression setting in which observations come from a single arm.  In much of the bandit literature \citep{chu2011contextual,bubeck2012regret,perchet2013multi,hu2022smooth}, effective algorithms rely on estimating unobserved rewards and deriving sharp non-asymptotic confidence bounds. To develop algorithms for single-index bandits, it is therefore essential to establish uniform estimation guarantees for regression under single-index feedback.

Let $f_0$ and $v_0$ denote the target link function and index, and consider the model in \eqref{eq:reg-form}:
\begin{equation}\label{eq:reg-form-offline}
	Y_t = g_0(X_t) + \xi_t = f_0(v_0^\top X_t) + \xi_t,
\end{equation}
where $f_0$, $v_0$, and $\xi_t$ satisfy the assumptions in Section~\ref{sec:set-up}. We first briefly review existing methods for single-index regression and discuss why they are insufficient for our bandit setting.

The classical optimality result for single-index regression was established by \cite{gaiffas2007optimal} using an aggregation approach. They showed that the minimax-optimal MISE rate $O\left(n^{-\frac{2\beta}{2\beta+1}}\right)$ can be achieved via empirical risk minimization. However, this aggregation-based guarantee does not extend to pointwise or uniform error bounds. A similar limitation applies to the temperature method used to derive PAC-Bayesian bounds in \cite{alquier2013sparse}.
To obtain pointwise and uniform estimation guarantees, a standard strategy is to consistently estimate the index $v_0$, followed by nonparametric regression using the resulting plug-in estimator. Two major approaches have been developed for estimating $v_0$: inverse regression and forward regression.

\textit{Inverse regression.}
Inverse regression exploits the idea that, when the distribution of $X_t$ exhibits certain structural properties (e.g., elliptical symmetry), the conditional distribution $X_t \mid Y_t$ reveals information about the subspace $\spa(v_0)$. Classical methods include sliced inverse regression (SIR) \citep{li1991sliced} and the sliced average variance estimator (SAVE) \citep{dennis2000save}.
A key assumption underlying inverse regression is the linear conditional mean (LCM) condition, which requires $\E[X_t \mid v_0^\top X_t]$ to be linear in $v_0^\top X_t$. Unfortunately, this condition does not hold in contextual bandit settings, where the selection of arms alters the conditional distribution of $X_t$ in a manner that breaks the geometric properties of $P_X$ necessary for the LCM condition.

\textit{Forward regression.} A prominent class is the Stein's Lemma-based methods \citep{brillinger2011generalized}, which exploit the identity $\E [Y_t  X_t]=\E \left[f^{\prime}\left(v_0^\top X_t\right)\right]  v_0$ valid for Gaussian $P_X$; these require Gaussian covariates and $\beta\ge 1$. In a similar spirit, gradient-based approaches such as ADE \citep{hristache2001direct} and rMAVE \citep{xia2006semi} estimate $\nabla g(x)$, which carries information about $v_0$, but rely on well-behaved gradients and typically require $\beta \ge 2$, or even $\beta\ge 3$ for higher-order methods.

In contextual bandit problems, these existing approaches are unsuitable for online estimation: they either rely on assumptions (e.g., the LCM condition or Gaussian covariates) that fail under bandit sampling, or they apply only to restricted smoothness regimes. Motivated by these limitations, we instead study a rank-based estimator, which imposes fewer structural requirements and is well-suited to our online setting.

\subsection{Maximum rank correlation and index plug-in}

Rank-type estimators are a class of M-estimators that learn latent index structures by maximizing or minimizing rank-based objective functions. A canonical example is Han's maximum rank correlation (MRC) estimator \citep{han1987non}, which recovers a single-index model under a monotone response.
For brevity, we write $Z_t = (X_t, Y_t)$ for the $t$-th observation, and denote the regression sample by $\cS = \left\{Z_t : t \in [n]\right\}$.
We now describe Han's MRC estimator \citep{han1987non}.
Define the kernel of the associated U-statistic as  $J(Z_1,Z_2, v)= \idc\left(Y_1>Y_2\right) \idc\left({X}_1^{\top} {v}>{X}_2^{\top} {v}\right)$. The rank-sum objective is then given by
\begin{equation}\label{eq:mrc}
	\begin{aligned}
		\Gamma_{\cS} (v) & = \frac{1}{\abs{\cS} (\abs{\cS} -1)} \sum_{Z_i\neq Z_j \in \cS} J(Z_1,Z_2, v) =   \frac{1}{n(n-1)} \sum_{ i\neq j } \idc\left(Y_i>Y_j\right) \idc\left({X}_i^{\top} {v}>{X}_j^{\top} {v}\right).
	\end{aligned}
\end{equation}
The MRC estimator is obtained by maximizing the empirical rank correlation:
\begin{equation}\label{eq:mrc-solve}
	\widehat{v}^{\sfH} =\underset{u\in\cV}{\operatorname{argmax}} ~\Gamma_{\cS}(u),
\end{equation}
where $\cV\subseteq\{u\in\R^d: u_1=1\}$ is a known compact set containing $v_0$. Without loss of generality, we take the response to be increasing in the index, as in~\eqref{eq:stoch-monotone}. The decreasing case can be handled analogously by replacing the kernel with $J'(Z_1,Z_2, v)= \idc\left(Y_1<Y_2\right) \idc\left({X}_1^{\top} {v}>{X}_2^{\top} {v}\right)$.

Our single-index regression procedure using MRC is summarized in Algorithm~\ref{alg:MRC-reg}, where the estimator $\widehat{v}^{\sfH}$ is used as input for local polynomial estimation (LPE).

\begin{algorithm}[h!]
	\caption{Single-Index Regression By MRC}
	\label{alg:MRC-reg}
	\begin{algorithmic}[1]
		\REQUIRE $\cS=\{(X_t,Y_t):t\in [n]\}$, where $n=\abs{\cS}$, smoothness level $\beta$, tuning parameter $C_{\sfH}$ .
		\STATE Split $\cS$ evenly into $\cS_1$, $\cS_2$
		\STATE Compute $	\widehat{v}^{\sfH} =\underset{u\in\cV}{\operatorname{argmax}} ~\Gamma_{\cS_1}(u) $ following \eqref{eq:mrc}, \eqref{eq:mrc-solve}
		\STATE For each $(X_t,Y_t)\in \cS_2$, project $X$ in onto direction $\widehat{v}^{\sfH}$ and get new  $\cS_2'=\{( (\widehat{v}^{\sfH})^{\top} X_t,Y_t): (X_t,Y_t)\in \cS_2\}$.
		\STATE Run local polynomial estimation on $\cS_2'$ with level $\lfloor \beta
		\rfloor$, bandwidth $h_n=\left(\frac{\log n}{n}\right)^{\frac{1}{2\beta+1}}\vee C_{\sfH}\left(\frac{d+\log^2 n  }{n} \right)^{\frac{1}{2}}$ and uniform kernel $\idc\{\abs{\cdot}\le 1  \}$, denote the regression function as $\widehat{f}(\cdot)$.
		\STATE Plug $	\widehat{v}^{\sfH} $ into $\widehat{f}(\cdot)$ and get  $\widehat{g}(x)=\widehat{f}( (\widehat{v}^{\sfH})^{\top} x )$.

		\STATE (Optionally) repeat 2-5 by switching $\cS_1$ and $\cS_2$ and get $\widehat{g}'$. Let $\widehat{g} \gets \frac{1}{2}(\widehat{g}+\widehat{g}') $ for cross-fitting.

		\STATE \textbf{return} Regression function $\widehat{g}(x)$
	\end{algorithmic}
\end{algorithm}

Instead of relying on assumptions that are difficult to maintain under bandit sampling, such as the LCM condition or Gaussian covariates, Algorithm~\ref{alg:MRC-reg} relies on the rank monotonicity of Assumption~\ref{asp:model}, a weak form of the monotone link commonly adopted in many reward models \citep{kakade2011efficient,han2017provable}. We stress that this is an exchange of assumptions rather than a weakening of them: the rank-based approach dispenses with smoothness requirements on the link ($\beta\ge1,2$ or $3$ for the forward-regression methods above) and with geometric requirements on $P_X$, but replaces them with rank monotonicity, with a bounded conditional density of $X_1$ given the remaining coordinates (Assumption~\ref{asp:sufficient-cond:1}) and with the nondegeneracy conditions of Assumption~\ref{asp:nondeg}. What makes the exchange attractive here is that the latter conditions are \emph{preserved} under the adaptive, arm-dependent sampling induced by an elimination policy, which is typically not the case for the LCM condition or for Gaussianity. Moreover, monotonic reward structures arise naturally in many applications: the Netflix artwork selection problem discussed earlier \citep{chandrashekar2017artwork,gur2022smoothness}, binary-outcome decision-making modeled with multinomial logit rewards \citep{mcfadden1987regression}, product and brand selection \citep{elshiewy2017multinomial}, and transportation choice modeling \citep{zhao2020prediction}.

\begin{Remark}[Role of monotonicity]\label{rmk:monotone} The rank monotonicity~\eqref{eq:stoch-monotone} enters our analysis only in the index estimation step, where it is of considerable help: it is precisely the condition under which the population rank correlation is maximized at the true index. Dropping it is not innocuous. Take $d=2$, $X\sim\cN(0,\bI_2)$ truncated to a ball, $v_0=(1,0)^\top$ and the non-monotone link $f_0(u)=u^2/c$ with additive symmetric noise, so that the response depends on the covariate only through $\abs{v_0^\top X}$. The sign flip $(X_1,X_2)\mapsto(-X_1,-X_2)$ leaves the design and the responses invariant while reversing the event $\{X_1^\top v>X_2^\top v\}$, whence $\Gamma(v)=\frac{1}{2}\PP(Y_1>Y_2)=\frac{1}{4}$ for \emph{every} $v$ with $v_1=1$: the criterion is flat and the MRC estimator cannot identify $v_0$ at any sample size. Without rank monotonicity one would need another criterion, or the forward-regression approaches above with their smoothness and design requirements. Whether sup-norm optimal index plug-in estimation is possible without it for $\beta<1$ appears to be open even offline.
\end{Remark}

\subsection{Theoretical properties of MRC and index plug-in}

Despite the long history of MRC estimators in the literature (see, for example, \cite{han1987non}, \cite{sherman1993limiting}, \cite{abrevaya2011rank}, \cite{han2017provable} and \cite{fan2020rank}), existing theories do not provide sharp non-asymptotic bounds for regression. The most recent analysis in \cite{fan2020rank} establishes only asymptotic convergence rates, which are insufficient for bandit regret control. Motivated by this gap, we develop in this section a new non-asymptotic theory for MRC that yields sharp finite-sample error bounds.

We write the covariate $X$ as $X = [X_1, \widetilde{X}^\top]^\top$, where $\widetilde{X} \in \R^{d-1}$ contains the last $d-1$ components of $X$. Correspondingly, we write $x = [x_1, \widetilde{x}^\top]^\top$, and let $\mu_0(\cdot)$ denote the marginal density of $v_0^\top X$. We introduce the following definition:
\begin{Definition}[$(c_0,r_0,r_1)$-regularity condition]\label{def:regular}
	A set $A \subseteq \R$ is called $(c_0,r_0,r_1)$\textit{-regular} if
	\begin{equation}
		\Leb(A \cap U(x,r)) \ge c_0 \cdot 2r,
		\qquad \forall\, r_1 < r \le r_0,\ \forall\, x \in A.
	\end{equation}
\end{Definition}
When $r_1 = 0$, we refer to this simply as the $(c_0,r_0)$-regularity condition.
Following \cite{audibert2007fast}, we say a distribution satisfies the \emph{strong density} condition if its support is a $(c_0,r_0)$-regular set for some constants $c_0, r_0$, and its density is bounded away from zero and infinity on that support.
The assumptions required for our analysis are provided in Assumption~\ref{asp:sufficient-cond}.
\begin{Assumption}[Conditions for estimation]\label{asp:sufficient-cond}
	We assume that:
	\begin{enumerate} [label=\ref*{asp:sufficient-cond}.\arabic*]
		\item \label{asp:sufficient-cond:1} For almost every $\widetilde{x}$, the conditional law of $X_1$ given $\widetilde{X}=\widetilde{x}$ has a density bounded by a constant $C_{\mu}$.

		\item  \label{asp:sufficient-cond:2}
		$v_0^\top X$ is with $(c_0,r_0)$-regular support $I_0:= v_0\circ\cX = \left\{v_0^\top x\mid x\in\cX\right\}$, and its marginal density $\mu_0(\cdot)$ satisfies $\lmu \le \mu_0(\cdot) \le \umu $ on its support $I_0$.
	\end{enumerate}
\end{Assumption}
	Assumption~\ref{asp:sufficient-cond} controls the regularity of the design, but it does not by itself guarantee that the index is identified at a parametric rate: the population rank correlation must also be \textit{nondegenerate} at $v_0$. In previous analyses of MRC \citep{sherman1993limiting,fan2020rank} this is expressed through the second-order matrices
	\begin{equation}\label{eq:well-cond-mat-main}
		\begin{aligned}
			\Delta_1(P) & = \E_P \left[\nabla_v\tau_P(Z,v_0)\nabla^\top_v\tau_P(Z,v_0) \right], \quad  \Delta_2(P)  & = -\E_P\left[\nabla^2_v\tau_P(Z,v_0)\right],
		\end{aligned}
	\end{equation}
	where $\tau_P(z,v) =\E_P J(z,Z,v)  + \E_P J(Z,z,v)$ for $(X,Y)\sim P$ and the derivatives $\nabla_v,\nabla^2_v$ are taken with respect to the index parameter $v$ (recall that $v_1=1$, so that $\nabla_v$ is a gradient in the free coordinates $\widetilde{v}\in\R^{d-1}$) and evaluated at $v=v_0$. Assuming directly that $\lambda_{\min}(\Delta_2(P))$ is bounded below is unsatisfactory, however, because the condition is stated in terms of second derivatives of a population $U$-statistic and is not transparent to verify. In particular it has to be checked not only for $P$ but also for the conditional laws induced by an elimination policy (Section~\ref{sec:bandits}) and for the instances used in the lower bound construction (Section~\ref{sec:opt-rates}). We therefore work instead with the following primitive conditions.
	Define the \textit{concordance function}
	\begin{equation}\label{eq:psi-def}
		\psi_0(w_1,w_2) := \tfrac{1}{2}\, \E\left[\sgn(Y_1-Y_2)\mid v_0^\top X_1 = w_1,\ v_0^\top X_2 = w_2 \right],
	\end{equation}
	where $(X_1,Y_1)$ and $(X_2,Y_2)$ are independent copies of $(X,Y)\sim P$. Under Assumption~\ref{asp:model}, $\psi_0$ is antisymmetric and $\psi_0(w_1,w_2)\ge 0$ whenever $w_1\ge w_2$.
	\begin{Assumption}[Nondegeneracy]\label{asp:nondeg} There exist constants $\cf,r_{\sff},\lmu_1,\ell_1, c_{\sfX}>0$ such that:
		\begin{enumerate} [label=\ref*{asp:nondeg}.\arabic*]
			\item \label{asp:nondeg:1} (Nondegenerate monotone response) $\psi_0(w_1,w_2)\ge \cf\,(w_1-w_2)$ for all $w_1\ge w_2$ in $I_0$ with $w_1-w_2\le r_{\sff}$.
			\item \label{asp:nondeg:2} (Conditional index density floor) For almost every $\widetilde{x}$ in the support of $\widetilde{X}$, the conditional law of $v_0^\top X$ given $\widetilde{X}=\widetilde{x}$ is absolutely continuous, and there is an interval $J_0\subseteq I_0$ with $\Leb(J_0)\ge \ell_1$, not depending on $\widetilde{x}$, on which its density is at least $\lmu_1$.
			\item \label{asp:nondeg:3} (Nondegenerate design) $\lambda_{\min}\left(\Cov(\widetilde{X})\right)\ge c_{\sfX}$.
		\end{enumerate}
	\end{Assumption}
	Assumption~\ref{asp:nondeg:1} holds with $\cf = \frac{1}{2}\inf_{I_0} f_0'$ for Bernoulli responses, since then $\psi_0(w_1,w_2)=\frac{1}{2}\left(f_0(w_1)-f_0(w_2)\right)$, and with $\cf \asymp \inf_{I_0} f_0'\cdot \lmu_u$ for additive noise regression whose noise difference has density at least $\lmu_u$ near the origin. Assumption~\ref{asp:nondeg:2} requires that, whatever the remaining coordinates may be, the index still takes every value in a fixed interval with non-vanishing density; it is not implied by a two-sided density bound for $X_1$ given $\widetilde{X}$. Assumption~\ref{asp:nondeg:3} asks that the design be nondegenerate in every direction. All three are conditions on the primitive distribution, and Section~A.4 of the Supplementary Material~\citep{suppSIB} shows that the conditional laws arising in the bandit problem inherit them with explicit constants.

	These conditions are exactly what is needed to make the population criterion $\Gamma(v)=\E J(Z_1,Z_2,v)$ nondegenerate. Provided the projections of $X$ are sub-Gaussian, $\norm{u^\top X}_{\psi_2}\le C_{\psi}$ for every $u\in\bbB_2(1)$, which holds with $C_{\psi}=2C_{\cX}$ under \eqref{eq:design-bounded}, the criterion separates globally and is curved locally at $v_0$:
	\begin{equation}\label{eq:curvature-main}
		\Gamma(v_0)-\Gamma(v) \ \ge \ c_{\Delta}\left( \norm{v-v_0}_2^2 \wedge r_{\sfc}^2 \right) \qquad \text{for every } v\in\cV,
	\end{equation}
	with curvature constant $c_{\Delta}= \frac{1}{4}\cf\lmu_1^2\ell_1 c_{\sfX}$ and a radius $r_{\sfc}$ proportional to $r_{\sff}\wedge\frac{\ell_1}{2}$; if the conditional density of $X_1$ given $(\widetilde{X},Y)$ is moreover three times boundedly differentiable, one also obtains $\lambda_{\min}(\Delta_2(P))\ge 2c_{\Delta}$ and $\lambda_{\max}(\Delta_2(P))\le C\lambda_{\max}(\Cov(\widetilde{X}))$. See Section~A.1 of the Supplementary Material~\citep{suppSIB} for more details. With it, the distribution class can be described purely through primitive conditions.
		\begin{Definition}[Distribution class]\label{def:dist-offline}
		Given constants $(C_{\mu},c_0,r_0,\lmu,\umu)$ in Assumption~\ref{asp:sufficient-cond} and $(\cf,r_{\sff},\lmu_1,\ell_1,c_{\sfX})$ in Assumption~\ref{asp:nondeg}, we denote by $\cP_0$ the collection of distributions $P$ for $Z=(X,Y)$ that satisfy Assumptions~\ref{asp:model}, \ref{asp:sufficient-cond} and \ref{asp:nondeg} with these constants. Every $P\in\cP_0$ then satisfies~\eqref{eq:curvature-main} with $c_{\Delta}$ determined by these constants.
	\end{Definition}
For fixed $d$, the constants of $\cP_0$, and hence $c_{\Delta}$, are constants of the problem.
The theoretical result of our Algorithm \ref{alg:MRC-reg} is then given by the following theorem:
	\begin{Theorem}[Estimation error]\label{thm:est} Under Assumptions \ref{asp:model}, \ref{asp:sufficient-cond:1} and \ref{asp:nondeg}, there exists a constant $C_{0}>0$ such that for every $\varepsilon\in(0,e^{-1}]$, with probability at least $1-\varepsilon$,
		\begin{equation}\label{eq:main-mcr-est}
			\norm{	\widehat{v}^{\sfH}  -v_0 }_2\le C_0\sqrt{\frac{d+\log^2(\frac{1}{\varepsilon})}{n}},
		\end{equation}
		provided that $n \ge C_{v,0}  \left(d+\log^2(\frac{1}{\varepsilon})\right)$, where $C_{v,0}=C_{\gap}/(c_{\Delta}r_{\sfc}^2)^2$.

		Moreover, if we further impose Assumption \ref{asp:sufficient-cond:2} and set $h_n=\left(\frac{\log n}{n}\right)^{\frac{1}{2\beta+1}}\vee C_{\sfH}\left(\frac{d+\log^2(\frac{1}{\varepsilon})  }{n} \right)^{\frac{1}{2}}$  for some large $C_{\sfH}$, then there exist constants $c_{\sfH}\gtrsim C_0/C_{\sfH}$ and $C_{v,\infty}\ge C_{v,0}$ such that, when \eqref{eq:main-mcr-est} holds and $n\ge C_{v,\infty}(d+\log^2(\frac1\varepsilon))$, domain $\widehat{v}^{\sfH} \circ X \subseteq U(I_0,c_{\sfH} h_n) $, and  with probability at least $1-\varepsilon$ we have the uniform bound:
		\begin{equation*}
            \begin{aligned}
                			\norm{ \widehat{f}_0- f_0 }_{\infty} & =  \sup_{U(I_0,c_{\sfH} h_n) }\abs{\widehat{f}_0 ( a)- f_0( a )}
                            \\
                            &\le C_{0,\infty} \left[ \sqrt{ \frac{\log(\frac{n}{\varepsilon})}{\log n}} \left( \frac{\log n }{n}\right)^{\frac{\beta}{2\beta+1}} +  \left(\frac{d+\log^2(\frac{1}{\varepsilon}) }{n} \right)^{\frac{\beta\wedge 1 }{2}}\right],
            \end{aligned}
		\end{equation*}
        some $ C_{0,\infty} >0$. The constants $C_{0}$, $C_{v,0}$, $C_{v,\infty}$ and $C_{0,\infty}$ depend only on $\beta$, $L$, $C_{\sfH}$, $C_{\cX}$, the radius of $\cV$ and the constants of Assumptions \ref{asp:model}, \ref{asp:sufficient-cond} and \ref{asp:nondeg}; in particular they are free of $n$, and they depend on $d$ only through these quantities.
	\end{Theorem}
	Notice that the cross-fitting in Algorithm \ref{alg:MRC-reg} is for robustness and better data efficiency, which will not change the order of error rate.

	Theorem~\ref{thm:est} is stated under the bounded design~\eqref{eq:design-bounded} for simplicity; both of its bounds in fact only require the projections of $X$ to be sub-Gaussian, at the price of a $\sqrt{\log n}$ factor that is absorbed in the polylogarithmic factors of our regret bounds. The latter version is used in Section~\ref{sec:increasing-dim}.

	An immediate consequence of Theorem \ref{thm:est} is that when $d\lesssim n^{\frac{1\wedge(2\beta -1)}{2\beta +1}}$, the following uniform bound holds for the final plug-in regression estimator $\widehat{g}_0(x)=\widehat{f}_0( (\widehat{v}^{\sfH})^{\top} x )$.
	\begin{Corollary}[Uniform bound]\label{coro:unif} Under the same condition as Theorem \ref{thm:est}, if $d \le  C n^{\frac{1\wedge(2\beta -1)}{2\beta +1}}$ for some constant $C>0$, then there exist constants $C_{0,\infty},\widetilde{C}_{1}>0$, again depending only on the constants of Assumptions \ref{asp:model}, \ref{asp:sufficient-cond} and \ref{asp:nondeg} and on $\beta,L$, such that with probability $1-n^{-10}$,
		\begin{equation}\label{eq:reward-est-uniform}
			\norm{\widehat{g}_0- g_0}_{\infty}= \sup_{x\in\cX}\abs{ \widehat{f}_0( (\widehat{v}^{\sfH})^{\top}x  )- f_0( v_0^\top x ) }  \le  C_{0,\infty} \left(\frac{ \log n }{n}\right)^{\frac{ \beta }{2\beta+1 } },
		\end{equation}
		provided that $n \ge \widetilde{C}_{1}  \left(d+\log^2 n  \right)$.
	\end{Corollary}
	Our Theorem \ref{thm:est} and Corollary \ref{coro:unif} show that  Algorithm \ref{alg:MRC-reg} achieves optimal single-index regression with sharp uniform non-asymptotic error bound, and is therefore suitable for our single-index bandits. Since $g_0$ takes values in $[0,1]$, clipping $\widehat{g}_0$ to $[0,1]$ can only reduce its error, and bounds it by $1$ on the exceptional event of probability $n^{-10}$; the clipped estimator $\widetilde{g}_0=(\widehat{g}_0\vee0)\wedge1$ therefore has the uniform $\ell_{\infty}$ risk
	\begin{equation*}
		\sup_{\cP_0,  \cH(\beta, L) }\E \norm{\widetilde{g}_0- g_0}_{\infty} \lesssim \left(\frac{ \log n }{n}\right)^{\frac{ \beta }{2\beta+1 } }.
	\end{equation*}
	This rate  matches the minimax optimal rate for single-index regression  \cite{gaiffas2007optimal,lepski2014adaptive}, which aligns with the minimax optimal rate for $1$-dimensional nonparametric regression as in \cite{stone1982optimal}, \cite{tsybakov2009introduction}.

	\section{Single-Index Bandits with Online Estimation}
	\label{sec:bandits}

	In this section, we present our main algorithm for solving single-index bandits, along with the corresponding assumptions that support our theoretical analysis.

	\subsection{Batched bandit algorithm without bin elimination}

	Theorem \ref{thm:est} provides a powerful tool to assess the uncertainty of offline single-index regression. It also gives us theoretical foundations to construct upper confidence bounds, a key tool to balance between exploration and exploitation in bandit problems, as long as we can perform online single-index regression following a similar fashion. Following this idea, we describe our batched online contextual bandit algorithm in Algorithm \ref{alg:sim-bandit}.

We stress at the outset that Algorithm~\ref{alg:sim-bandit} is an \textit{epoch-based batched} policy rather than a fully online procedure: the estimates $\widehat{g}_{k,m}$ and the active arm sets $\cK_m(\cdot)$ are updated only at the end of each epoch, and within an epoch the decision rule is frozen. This is not merely a convenience: it is what makes the data collected within an epoch i.i.d.\ conditionally on the history, and hence makes the MRC analysis of Section~\ref{sec:regression} applicable to adaptively collected data. Nor is the algorithm of the usual fixed-batch-budget form: the number of epochs is not specified a priori but determined by $\{\epsilon_m\}$, the epoch lengths $\{n_m\}$ and the horizon $n$, yielding $M=O(\log n)$ epochs.
	\begin{algorithm}[h!]
		\caption{Batched Contextual Bandit with Single-Index Rewards}
		\label{alg:sim-bandit}
		\begin{algorithmic}
			\REQUIRE Epoch scheme $\{\cT_m\}_{m=1}^{M}$, confidence bound $\{\epsilon_m\}_{m=0}^{M-1}$  smoothness level $\beta$, $\widehat{g}_{k,0}=0$, active arm set function $\cK_{0}(x)=[K]$
			\FOR{$m \in \{1,\dots,M\}$}

			\FOR{$t\in\cT_m$}
			\STATE Receive covariate $X_t$, and compute active arm set $\cK_{m-1}(X_t)$
			\IF{$\abs{\cK_{m-1}(X_t)}\ge 2$ }
			\STATE  Sample arm $\pi_t\in\cK_{m-1}(X_t)$ by uniform sampling
			\ELSE
			\STATE  Choose the only arm in $\cK_{m-1}(X_t)$
			\ENDIF
			\ENDFOR

			\STATE Log data for each arm $k$: $\cS_{k,m}= \{(X_t,Y_t)\mid t\in \cT_m,\pi_t=k\}$.
			\STATE For each arm $k$, update estimation $\widehat{g}_{k,m}$ by running Algorithm \ref{alg:MRC-reg} on $\cS_{k,m}$ with smoothness $\beta$.
			\STATE Pre-select arm by:
			\begin{equation*}
				\cK_{m-0.5}(x)=\left\{k\in \cK_{m-1}(x)\mid \widehat{g}_{(1),m-1}^{}(x, \cK_{m-1}(x))- \widehat{g}_{k,m-1}(x )\le \frac{1}{2}\epsilon_{m-1} \right\}
			\end{equation*}
			\STATE Update possible arm set
			\begin{equation*}
				\cK_{m}(x)=\left\{k\in\cK_{m-0.5}(x)\mid \widehat{g}_{(1),m}^{}(x, \cK_{m-0.5}(x))- \widehat{g}_{k,m}(x )\le \epsilon_{m}  \right\}
			\end{equation*}

			\ENDFOR
			\STATE \textbf{return} $\{\pi_t\}_{t=1}^n$
		\end{algorithmic}
	\end{algorithm}

	Our Algorithm \ref{alg:sim-bandit} features a batched structure and an arm pre-selection approach before we update estimation for the next epoch.
	We will show later that this pre-selection approach plays a crucial role in learning rewards near the margin region: the pre-selection step, informed by data from previous epochs, discards those arms whose reward would otherwise have to be estimated on irregular regions of the context space. Unlike most existing methods \citep{perchet2013multi,hu2022smooth,gur2022smoothness,cai2024transfer}, our algorithm does not rely on bin partition and elimination. Instead, it implicitly defines decision regions through $\cK_{m}(x)$, thereby avoiding the $2^d$ factor in the regret upper bound that arises from explicit bin partition.

	We explain some notation and settings used in Algorithm \ref{alg:sim-bandit}. For any arm set $\cK$, we define $\widehat{g}_{(1),m}^{}(x, \cK)= \max\limits_{k\in \cK } \left\{\widehat{g}_{k,m}(x) \right\}$, and $\widehat{g}_{(2),m}(x,\cK)$ is defined similarly following $g_{(2)}(x)$ in \eqref{eq:margin-def}. Here, we select the batch and confidence bound with $\epsilon_m= c_{\epsilon} 2^{-m}$ with a small $c_{\epsilon}$, and let
	\begin{equation*}
		n_m=\left\lceil C_{\cT} \left(\frac{d+\log^2 n }{\epsilon_m^{\frac{2}{1\wedge \beta}}  }+  \left(\frac{ \log n}{\epsilon_m^2}\right)^{\frac{2\beta+1}{2\beta}} \right) \right\rceil
	\end{equation*}
	be the length of each epoch for some large $C_{\cT}$, with $n_0=0$.
	Let $M=\inf\{m\mid \sum_{i=1}^m n_{i} \ge n\}$ be the total number of epochs. We then define  $\cT_m=\{ t \mid \sum_{i=1}^{m-1} n_{i}+1 \leq t \leq \min \left\{\sum_{i=1}^m n_i,\ n\right\}\}$, for $m\in [M]$, and the number of accumulated observations after epoch  $\cT_m$ as $S_m= \sum_{i=1}^m n_i$. Since $\epsilon_m$ decays geometrically, we have $n_m \asymp 2^{\frac{2m}{1\wedge \beta}}(d+\log^2 n)\vee 2^{\frac{m(2\beta+1)}{\beta}}\log^{\frac{2\beta+1}{2\beta}} n$ and therefore $M= O(\log n)$.

	\subsection{Assumptions for single-index bandits}\label{sec:bandit-asp}

	To move from offline to online estimation and bandit regret control, more assumptions on the optimal decision regions are required so that the learning problem itself is not pathological.  We denote
	$$\cQ_k=\{x\in \cX\mid  g_{k}(x)= g_{(1)}(x) \},  \quad I_k:= v_k\circ \cQ_{k}  = \{v_k^\top x\mid x\in \cQ_k\} $$
	as the optimal decision region and it projection on index $v_k$ for each arm $k\in[K]$. The conditions below are imposed on $\cQ_k$ and $I_k$ \textit{for every arm} $k\in[K]$ simultaneously:
	\vspace{2pt}

	\begin{Assumption}[Learnable margin area]\label{asp:decision-region}For every arm $k\in[K]$, we make the following assumptions:
		\begin{enumerate} [label=\ref*{asp:decision-region}.\arabic*]
			\item \label{asp:decision-region:1}  For each $k$ and any $x_1,\widetilde{x},y$ in the support of $(X,Y^{(k)})$, $ \PP(X\in \cQ_k\mid \widetilde{X}=\widetilde{x}, Y^{(k)}=y) \ge \lp_K $.

			\item \label{asp:decision-region:2} Define  margin area for arm $k$ as  $\cQ^{\sfm}_k(\delta):=\{x\mid 0<g_{(1)}(x)-g_{k}(x) \leq \delta\}$, with its projection on $v_k$ as $I_k'(\delta) := v_k\circ \cQ^{\sfm}_k(\delta)$. There exist constants $\delta_{1}, c_{\sfm},\lmu_{\sfm}>0$ such that, for any $n^{-\frac{\beta}{2\beta+1}}\le \delta\le \delta_1$, the set $I_k\cup I_k'(\delta)$ is $\left(c_0,c_{\sfm}\delta^{1/\beta}, c_{\sfm} n^{-\frac{1}{2\beta+1}} \right) $-regular, and conditional on $X\in \cQ_{k}\cup \cQ^{\sfm}_k(\delta) $, the variable $v_k^\top X $ has density that satisfies $\mu_k(\cdot\mid \cQ_{k}\cup \cQ^{\sfm}_k(\delta)) \ge\lmu_{\sfm} $ on its support.

			\item \label{asp:decision-region:3} (Nondegeneracy on the optimal region) The conditional law $P^{(k)}\mid \cQ_k$ of $(X,Y^{(k)})$ given $X\in\cQ_k$ satisfies Assumption~\ref{asp:nondeg} with index $v_k$ and with the same constants $(\cf,r_{\sff},\lmu_1,\ell_1,c_{\sfX})$. In particular $\lambda_{\min}\left(\Cov\left(\widetilde{X}\mid X\in \cQ_k\right)\right)\ge c_{\sfX}$.
		\end{enumerate}
	\end{Assumption}

	Assumption \ref{asp:decision-region:1} requires that the (conditional) probability of $\{X \in \cQ_k\}$ is bounded away from zero. This rules out the pathological cases where, conditional on optimal decisions $\{X \in \cQ_k\}$, the distribution of $X_1| (\widetilde{X},Y^{(k)})$ is degenerate.  Due to the linear structure of single-index reward, this condition is relatively easy to satisfy. For example, when $Y^{(k)}$ is Bernoulli and $\min\{f_k, 1-f_k\}\ge c$, with  $v_1 = \cdots = v_K$, a sufficient condition is that  $\PP(X_1+ \widetilde{v}^\top \widetilde{X} \in I_k\mid \widetilde{X})\ge c$,  which holds, for instance, when $I_k$ is non-vanishing and $X_1$  has a density bounded away from zero conditional on $ \widetilde{X} $.

	Assumption \ref{asp:decision-region:2} ensures that the margin area is regular and therefore learnable for LPE; it holds when $I_k$ is regular and $I_k'(\delta)$ is either connected to $I_k$ or itself regular with length of order $\delta^{1/\beta}$.

Assumption~\ref{asp:decision-region:3} is the bandit counterpart of Assumption~\ref{asp:nondeg}: it requires the rank-correlation criterion to remain nondegenerate after restricting to the contexts on which arm $k$ is optimal, which is what allows $v_k$ to be estimated from the data actually collected for arm $k$. Being stated in terms of primitive quantities, it can be checked directly on the conditional law: Assumption~\ref{asp:nondeg:1} is inherited automatically, since $\psi_0$ depends only on the conditional law of $Y^{(k)}$ given the index, and Assumption~\ref{asp:nondeg:3} reduces to $\cQ_k$ having non-empty interior with the design density bounded below there. Section~E of the Supplementary Material~\citep{suppSIB} proves the latter reduction and verifies Assumptions~\ref{asp:sufficient-cond}, \ref{asp:nondeg}, \ref{asp:decision-region} and \ref{asp:margin} in full for two canonical families of instances, with a shared and with arm-specific indices.

Assumption~\ref{asp:decision-region:1} requires every arm to be optimal on a conditionally non-negligible set, so that instances in which some arm is never optimal are excluded from $\cP$, although they may well be easier.

	We require Assumption \ref{asp:decision-region:2}  as establishing learnability guarantees in single-index bandits is more challenging than in general nonparametric bandits. In the latter case, the optimal decision region is a $d$-dimensional subset of $\cX \subseteq \R^d$, which can be discretized using bin partition. When $P_X$ has strong density in $\R^d$, the discretized cubes are learnable at the expense of a $2^{d}$ factor in the regret upper bound; see, e.g., \cite{perchet2013multi,hu2022smooth}. By contrast, for single-index bandits, the optimal decision region lies in an \textit{unknown} subspace of $\R^d$, rendering gridding schemes inefficient. Furthermore, since $\{v_k\}$ are unknown, it is impossible to identify irregular regions and thereby perform the screening procedure proposed in \cite{hu2022smooth}.

	To reach small regret for nonparametric bandit, the following standard margin assumptions are required on the covariate distribution $P_X$ \citep{perchet2013multi,hu2022smooth,cai2024transfer}:
	\begin{Assumption}[Margin condition]\label{asp:margin}
		There exist constants $\alpha \geq 0, \delta_0>0,C_\alpha>0$ such that the reward functions $\left\{g_k\right\}_{k=1}^K$ and marginal distribution $P_X$ satisfy
		$$
		\PP_{X}\left(0<g_{(1)}(X)-g_{(2)}(X) \leq \delta\right) \leq C_\alpha \delta^\alpha, \quad \forall \delta\in  (0, \delta_0].
		$$
	\end{Assumption}

	\begin{Definition}[Bandit distribution class]\label{def:dist-bandit}  We denote by $\cP_{\alpha,\beta}$ (or simply $\cP$ when there is no ambiguity) the class of bandit distributions such that, for each arm $(X, Y^{(k)}) \sim P^{(k)}$, Assumptions~\ref{asp:model}, \ref{asp:sufficient-cond} and~\ref{asp:nondeg} hold with a common set of constants, and the joint distribution satisfies Assumptions~\ref{asp:decision-region} and~\ref{asp:margin}. By the curvature bound~\eqref{eq:curvature-main}, every arm distribution $P^{(k)}$ and every conditional law $P^{(k)}\mid\cQ_k$ then obeys the curvature bound~\eqref{eq:curvature-main} with the same constant $c_{\Delta}$.
	\end{Definition}

	\section{Optimal Rates for the Regret}
	\label{sec:opt-rates}

	In this section, we study the performance of Algorithm \ref{alg:sim-bandit}, together with the minimax rate of single-index bandits to show that Algorithm \ref{alg:sim-bandit} is simple yet optimal.

	\subsection{Regret upper bound}

	We first give a theoretical analysis on the regret of Algorithm \ref{alg:sim-bandit} under the case when $\beta$ is known. We will discuss the adaptation to unknown smoothness level in Section \ref{sec:adaptation}.

	The tuning parameters are important for selecting epoch scheme and confidence bound. We select our tuning parameters such that $C_{\cT}$ exceeds a constant determined by the class $\cP$, and  $c_{\epsilon}\le 0.6 \left(\delta_0 \wedge\delta_1 \wedge 1\right)$, where $\delta_1$ is from Assumption \ref{asp:decision-region} and $\delta_0$ from Assumption \ref{asp:margin}.

Since the theory requires $c_{\epsilon}$ small and $C_{\cT}$ large relative to constants unknown to the decision-maker, Section~D.1 of the Supplementary Material~\citep{suppSIB} describes how these, and the bandwidth inflation factor $C_{\sfH}$, can be calibrated in practice.

	A regret upper bound of  Algorithm \ref{alg:sim-bandit}, denoted by $\widehat{\pi}$, is described in Theorem \ref{thm:minimax-ub}:
	\begin{Theorem}[Upper bound]\label{thm:minimax-ub} Suppose $d\lesssim n^{\frac{1\wedge(2\beta-1)}{2\beta+1}}$. For distributions in $\cP$ given in Definition \ref{def:dist-bandit}, and link functions in $\cH(\beta,L)$, our Algorithm \ref{alg:sim-bandit} achieves the regret
		$$
		\sup_{\cP,\cH(\beta,L)}	\Reg(\widehat{\pi})   \leq C (\log n)^{ \iota_1(\alpha,\beta) } \cdot \left[ n \left(\frac{ 1 }{n}\right)^{\frac{ \beta(1+\alpha) }{2\beta+1 } }   \vee 1\right] + C(d+\log^2 n)\log n,
		$$
		where $C>0$ depends only on the constants defining $\cP$ and on $\beta,L,\alpha$, so that it is free of $n$ and depends on $d$ only through those constants; its dependence on the number of arms is polynomial in $K/\lp_K$, with $\lp_K$ from Assumption~\ref{asp:decision-region:1} (Section~D.2 of the Supplementary Material~\citep{suppSIB}).

	\end{Theorem}

	Here $\iota_1$ only depends on $\alpha,\beta$, and $\iota_1(\alpha,\beta)=\frac{1}{2\beta}+\frac{\beta(1+\alpha)}{2\beta+1} $ when $\beta+1-\alpha\beta > 0$; for $\beta+1-\alpha\beta \le 0$ the bound holds with $\iota_1 = \frac{1}{2\beta}+\frac{\beta(1+\alpha)}{2\beta+1}+1$, and in this regime the bracket equals $1$. The last term is the cost of the first epochs, in which every arm is explored. When $\alpha(\beta\wedge1)\le1$, it is absorbed into the first term under the dimension condition of the theorem, and the bound reads $C(\log n)^{\iota_1}\, n^{1-\frac{(1+\alpha)\beta}{2\beta+1}}$. This is the regime studied in the nonparametric bandit literature; see, e.g., \cite{rigollet2010nonparametric,gur2022smoothness}.%

	\begin{Remark} Compared with general nonparametric bandits, whose minimax regret is of order $n^{1 - \frac{\beta(1+\alpha)}{2\beta + d}}$, the leading term of Theorem~\ref{thm:minimax-ub} carries no $d$ at all; the leading constant is a constant of the class $\cP$, which can also be free of $d$ (see Assumption~\ref{asp:normalization}). It also avoids the $(\log^{1/\beta} n \vee 2)^{d}$ factor that explicit bin partitions produce \citep{perchet2013multi,hu2022smooth,gur2022smoothness,cai2024transfer}, because the decision regions are defined implicitly through $\cK_{m}(x)$.

		Although we assume $K=O(1)$ throughout, the dependence of the bound on $K$ is explicit in our proof and polynomial, as displayed in the theorem, so that the guarantee remains meaningful for slowly growing $K$.
	\end{Remark}

	We now explain why the pre-selection procedure is important for regret control. Let $\cD_{k,m-1}=\left\{x\mid k\in \cK_{m-1}(x)\right\}$ be the region on which arm $k$ may be pulled in epoch $m$. Controlling the regret requires estimating $g_k$ on $\cD_{k,m-1}$ with increasing accuracy, but $\cD_{k,m-1}$ itself need not be regular for the LPE step of Algorithm \ref{alg:MRC-reg}, which would leave the estimation uncontrolled.
	To overcome this, we make use of the regular margin condition in Assumption \ref{asp:decision-region:2}: this condition ensures that for each $\cD_{k,m-1}$, there exists a smaller margin area $\cQ_{k}\cup\cQ_{k}^{\sfm}(\frac{3\epsilon_{m-1}}{4})\subseteq \cD_{k,m-1}$ that is regular, and that accurate single-index regression guarantees $\cD_{k,m-1}\subseteq \cQ_{k}\cup\cQ_{k}^{\sfm}(\frac{5\epsilon_{m-1}}{4} )$ for a slightly larger margin area. We pre-select arms using $\widehat{g}_{k,m-1}$, which is accurate in a broader margin $\cQ_{k}\cup\cQ_{k}^{\sfm}(\frac{5\epsilon_{m-1}}{4} )$ that contains $\cD_{k,m-1}$. After pre-selection, the region requiring finer decision-making is confined within $\cQ_{k}\cup\cQ_{k}^{\sfm}(\frac{3\epsilon_{m-1}}{4} )$ as long as the estimation error is sufficiently smaller than $\epsilon_{m}$.

	\subsection{Regret lower bound}

	Although Theorem \ref{thm:minimax-ub} presents an exciting regret upper bound that can achieve sufficient dimension reduction and matches the minimax optimal regret for one-dimensional nonparametric bandits, it is still questionable on whether the upper bound is optimal for the distribution class $\cP, \cH(\beta,L)$. In this section, we construct a lower bound showing that, when $\alpha\beta\le1$ or $\alpha=1\le\beta$, the regret upper bound given in Theorem \ref{thm:minimax-ub} is minimax optimal up to logarithmic factors. We describe it in Theorem \ref{thm:minimax-lb}:
	\begin{Theorem}[Regret lower bound]\label{thm:minimax-lb} Suppose that $\alpha\beta\le 1$, or that $\alpha=1\le\beta$.
		For any admissible policy $\pi$, one has the minimax lower bound
		\begin{equation}\label{eq:minimax-lb-main}
			\inf _\pi \sup _{\cP,\cH(\beta,L)} \Reg(\pi)\geq c\left( n \left(\frac{ 1 }{n}\right)^{\frac{ \beta(1+\alpha) }{2\beta+1 } }   \vee 1 \right)
		\end{equation}
		where $c>0$ is a constant independent of $n,d$.
	\end{Theorem}

	This result shows that, in this range, the derived regret rate is tight up to logarithmic factors. The proof is based on constructing instances that are difficult to distinguish. For $\alpha\beta\le1$, we perturb the reward difference by $W\asymp\delta^{\alpha\beta-1}$ bumps of height $\delta^{\beta}$, each of which can be switched independently and reverses the optimal arm on a set of contexts of probability of order $\delta$, and combine them through Assouad's lemma with $\delta\asymp n^{-\frac{1}{2\beta+1}}$; for $\alpha=1\le\beta$, a single perturbation that shifts a transversal crossing suffices. In addition to an instance construction that is distinct from those in \cite{hu2022smooth,cai2024transfer}, we need every instance to lie in $\cP$, so that in particular its indices remain identifiable at a parametric rate. For $\alpha\beta\le1$ we therefore perturb the means of the rewards without changing their ranks, so that all instances share one concordance function $\psi_0$; the resulting mean links need not be monotone, which Assumption~\ref{asp:model} allows. Section~A.5 of the Supplementary Material~\citep{suppSIB} verifies all the conditions of $\cP$ for every instance.

	Together with Theorem~\ref{thm:minimax-ub}, Theorem~\ref{thm:minimax-lb} shows that, over $\cP$, in its range of $(\alpha,\beta)$ and under the dimension condition of Theorem~\ref{thm:minimax-ub}, the fundamental difficulty of single-index bandits is the same as that of one-dimensional nonparametric bandits. The rationale is that, for responses that are monotone in the index, the indices of all arms can be learned at the parametric rate $\sqrt{d/n}$ without knowledge of the links (Theorem~\ref{thm:est}); after the index plug-in, the problem effectively reduces to a one-dimensional nonparametric bandit. When $d$ grows faster, learning the indices becomes the bottleneck, and Section~\ref{sec:increasing-dim} shows that it adds the parametric term $n(d/n)^{\frac{1+\alpha}{2}}$.

	\begin{Remark}[When is the oracle policy static?]
		Under a strong density condition on $P_X$ the threshold at which a single arm becomes optimal for almost every context is known to be $\alpha(\beta\wedge 1)>d$ for $K$-armed nonparametric bandits \citep{perchet2013multi}, and can be refined to $\alpha(\beta\wedge 1)>1$ for $K=2$ \citep{rigollet2010nonparametric,hu2022smooth,gur2022smoothness}. However, that condition may fail for $P_X$ in our setting while regularity holds for the marginals of $v_k^\top X$. Relaxing it to $\PP_X(\cX\cap\bbB_2(x,\varepsilon))\ge c_{\gamma}\varepsilon^{d+\gamma}$ on a convex support, Section~B of the Supplementary Material~\citep{suppSIB} shows that for every $K\ge 2$ the oracle is static when $\alpha(\beta\wedge 1)>1+\gamma$, a threshold attained when $\beta\ge1$, and that it improves to $\alpha(\beta\wedge 1)>1$ when the arms share a common index: for $\beta\ge1$ it deteriorates by exactly the amount $\gamma$ by which the density condition is relaxed.
	\end{Remark}

	\section{Adaptation to Unknown Smoothness}
	\label{sec:adaptation}

	This section discusses the adaptation to the unknown smoothness level $\beta$ in single-index bandits. Our Algorithm \ref{alg:sim-bandit} and theoretical upper bound established in Section \ref{sec:opt-rates} necessarily require knowing the exact smoothness level $\beta$. However, in many applications, $\beta$ is unavailable to decision makers. In this section, we discuss the adaptation to the unknown smoothness level when only a range of $\beta$, $[\betamin,\betamax]$ is given such that $\beta\in [\betamin,\betamax]$. We will show that under the standard self-similarity condition, achieving adaptation to unknown $\beta$ is possible by an adaptive algorithm whose regret still attains the optimal rate up to logarithmic factors, for the ranges of $\beta$ and $d$ stated in Theorem~\ref{thm:adaptive-regret-ub}.

	\subsection{Impossibility of adaptation}

	It is well established in the nonparametric multi-armed bandit literature that adaptation to an unknown smoothness level without incurring additional cost is impossible \citep{locatelli2018adaptivity, gur2022smoothness, cai2022stochastic}. Without further restrictions, no algorithm can simultaneously achieve the minimax-optimal rate across a family of nonparametric bandit problems with varying smoothness parameters. This impossibility result naturally extends to single-index bandits by considering the special case $v_k^\top X = X_1$, which reduces the problem to a one-dimensional nonparametric bandit. Adapting the argument of \cite{locatelli2018adaptivity,gur2022smoothness}, we show that cost-free adaptation is also impossible within the class $\cP, \cH(\beta, L)$.
	As discussed in \cite{gur2022smoothness}, this failure of adaptation in nonparametric bandits arises partly from the incomplete feedback inherent to the bandit setting. This contrasts with nonparametric regression \citep{hengartner1995finite, cai2013adaptive}, where adaptation can sometimes be achieved under shape constraints on the function class.

	\begin{Theorem}\label{thm:imp-adap} Let $0<\beta_1\le 1$, $\beta_1<\beta_2$ and $0\le\alpha\le 1$, and write $\opt_{\Reg}(n,\alpha,\beta):= n\cdot (\frac{1}{n})^{\frac{\beta(\alpha+1) }{2\beta+1}}$. If an admissible policy $\pi$ satisfies $\sup_{\cP_{\alpha,\beta_2},\cH(\beta_2,L)}\Reg(\pi)\le C\,\opt_{\Reg}(n,\alpha,\beta_2)$ for every $n$, then there is a constant $c>0$, independent of $n$ and $d$, such that
		$$
		\sup _{\cP_{\alpha,\beta_1},\cH(\beta_1,L)}\Reg(\pi )  \ \geq\ \frac{c\,n}{\opt_{\Reg}(n,\alpha,\beta_2)^{1+\beta_1}}\ =\ c \cdot \opt_{\Reg}(n,\alpha,\beta_1) \cdot n^{ \frac{(1+\alpha) \beta_1}{2 \beta_1+1}-\frac{(1+\beta_1)(1+\beta_2-\alpha\beta_2)}{2\beta_2+1} }.
		$$
	\end{Theorem}

	The exponent of the last factor increases with $\beta_2$, and it is positive for large $\beta_2$ exactly when $2(1+\alpha)\beta_1>(1-\alpha)(1+\beta_1)(2\beta_1+1)$, which for $\beta_1=1$ means $\alpha>\frac12$. For $\alpha=1$ it equals $\frac{2\beta_1}{2\beta_1+1}-\frac{1+\beta_1}{2\beta_2+1}$, which is positive for every $\beta_2>\beta_1=1$. For instance, when $\beta_1=1$, $\beta_2 =2$ and $\alpha =1$, a policy that is rate-optimal over $\cP_{\alpha,\beta_2}$ has worst-case regret over $\cP_{\alpha,\beta_1}$ exceeding $\opt_{\Reg}(n,\alpha,\beta_1)$ by a factor of order $n^{4/15}$. In this regime, a rate-optimal single-index bandit algorithm designed for $\cH(\beta_2, L)$ cannot be minimax-optimal for $\cH(\beta_1, L)$ without additional restrictions.

	\subsection{Self-similarity condition}

	Given the impossibility of adaptation to unknown $\beta$ under general settings, we seek to explore the adaptive algorithms when more assumptions are imposed. As suggested by \cite{gur2022smoothness, cai2022stochastic,cai2024transfer}, when the reward functions are ``self-similar'', it is possible to extend the smoothness estimation technique to nonparametric bandits and achieve adaptation on $\beta$. We now show that adaptation to $\beta$  is possible for single-index bandits under self-similar rewards. We describe the self-similarity condition introduced in \cite{picard2000adaptive}, \cite{gine2010confidence} by first defining the $\ell_2$-projection to polynomials:
	\begin{Definition}
		Let $f(\cdot)$ be any function and let $l,p\ge 0$ be integers. For a link $f_k$ and an interval $U$, we write $\Proj_h^p f_k(\cdot; U)$ for the $\ell_2$-projection of $f_k$ onto the space of polynomials of degree at most $p$ over $U$, measured with respect to the law of the index $v_k^\top X$. Concretely, for each $x\in U$ we define
		$$
		\Proj_h^p f_k(x ; U):=g(x), \quad \text { s.t. } \quad g=\arg \min _{q \in \operatorname{Poly}(p)} \int_U|f_k(u)-q(u)|^2 K\left(\frac{x-u}{h}\right) p_k(u \mid U) d u,
		$$
		where $p_k(\cdot\mid U)$ is the density of $v_k^\top X$ conditionally on $v_k^\top X\in U$,
		where the kernel is $K(\cdot)=\idc\left\{|\cdot| \leq 1\right\}$ with bandwidth $h$, and $\operatorname{Poly}(p)$ represents the collection of polynomials of degree no greater than $p$.
	\end{Definition}
	To facilitate the definition of self-similarity, we denote the lattice in $\R$ as $\cB_{2,l}=\{B_i\mid B_i =[\frac{i-1}{2^l},\frac{i}{2^l}],i\in\bbZ\}$. The self-similar functions are given in Definition \ref{def:ss}.

	\begin{Definition}[Self-similar functions]\label{def:ss}
		Given some $\beta \in[\betamin, \betamax]$, and finite constants $l_0 \geq 0$ and $b>0$, define the class of self-similar sets of function, $\cH^{s s}\left(\beta,L, b, l_0\right)$, to be the collection of the sets of functions $\left\{f_k\right\}_{k \in [K]}$ such that $f_k \in \mathcal{H}(\beta,L), k \in [K]$, and for which for any integers $l \geq l_0$, and $\lfloor\beta\rfloor \leq p \leq\lfloor\betamax\rfloor$, one has

		$$
		\max _{k\in [K] } \max _{ B \in \mathcal{B}_{2,l},\, B\subseteq v_k\circ\cX}  \sup _{x \in B}\left| \Proj_{2^{-l}}^p f_k(x ; B)-f_k(x)\right| \geq b \cdot 2^{-l \beta} ,
		$$
		where $v_k\circ\cX$ is the support of $v_k^\top X$.
	\end{Definition}

	The self-similarity condition complements the Hölder smoothness in the sense that the function class is hard to estimate with a bias lower bounded by $h^{\beta}$. Our first result for the self-similar reward class is that the self-similarity condition will not reduce the complexity of the problem.

	\begin{Theorem}[Regret lower bound under self-similarity]\label{thm:minimax-lb-ss} Under the assumptions of Theorem \ref{thm:minimax-lb}, for suitable constants $b$ and $l_0$, the minimax lower bound satisfies
		\begin{equation}\label{eq:minimax-lb-main-2-ss}
			\inf _\pi \sup _{\cP,\cH^{ss}\left(\beta,L, b, l_0\right)} \Reg(\pi)\geq c\left( n \left(\frac{ 1 }{n}\right)^{\frac{ \beta(1+\alpha) }{2\beta+1 } }   \vee 1 \right)
		\end{equation}
	\end{Theorem}
	The constants $b$ and $l_0$ do not depend on $n$: the hard instances of Theorem~\ref{thm:minimax-lb} are made self-similar by adding to every link a common term with an exact scaling of order $\beta$, placed away from the perturbations.

	\subsection{Adaptation under self-similarity condition}\label{sec:adapt-ss}

	We study the problem of adaptation to an unknown smoothness parameter $\beta$, assuming that only an interval $[\betamin, \betamax]$ is known such that $\beta \in [\betamin, \betamax]$. We assume $\frac{1}{2}<\betamin < 1 < \betamax$. Our idea is to first estimate $\beta$ by adapting Lepski’s method \citep{lepski1997optimal}, which identifies the optimal bandwidth (corresponding to the true smoothness) by comparing the estimation bias and stochastic error across different bandwidths during the exploration phase. We then use an undersmoothed estimator based on the estimated $\beta$ as an input to Algorithm~\ref{alg:sim-bandit} for regret control. The MRC and index plug-in components are employed for dimension reduction.

	We introduce several notations and settings. We denote
	\begin{equation}\label{eq:bandwidth-lepski}
		l_1=l:=\left\lceil\frac{\betamin \log _2 n}{(2 \betamax+1)^2}\right\rceil, \quad l_2=l+\left\lceil \frac{1}{\betamin} \log_2 \log n\right\rceil.
	\end{equation}
	Our approach focuses on comparing $\max_{B\in \cB_{2,l},k\in [K] }\abs{\widehat{f}_k(z;B, 2^{-l_1})-\widehat{f}_k(z; 2^{-l_2})}$, where the first estimate is fitted on the data in the cell $B$ and the second on all the data. We consider first estimating the smoothness parameter with $N_0$ data, and then plug the estimate $\beta_{\est}$ into Algorithm \ref{alg:sim-bandit}. To this end, we define the grid for comparison as $\cM=\left\{k2^{-l_3}, k\in\bbZ  \right\}$, where $l_3 = \left(1+\frac{\betamax}{\betamin}\right) l_1 + \frac{1}{\betamin}\log_2 \log n$, and $\cM(B)= \cM\cap B$ for any $B\in\cB_{2,l_1}$
	We propose the resulting estimator of the smoothness, which returns $\beta_{\est}=\frac{1}{l_1}\log_2\frac{1}{b_{\max}}-C_{\ell}\frac{\log_2\log n}{\log_2 n}$ from the largest discrepancy $b_{\max}$ over the grid, truncated to $[\betamin,\betamax]$; the pseudocode is Algorithm~4 in Section~A.8 of the Supplementary Material~\citep{suppSIB}.

	Its performance is guaranteed by the following proposition.
	\begin{Proposition}\label{prop:smooth-est} With  the choice $N_0= \lceil C_{\gap}(d+\log^2 n) n^{\frac{2\betamin(\betamax+1)}{(2\betamax+1)^2}} \rceil$ for some large $C_{\gap}>0$, we have, for $n\ge n_0$, that
		$$
		\mathbb{P}\left(\beta_{\est}  \in\left[\beta-\frac{C_{\est}(2 \betamax+1)^2 \log _2 \log n}{\betamin \log _2 n}, \beta\right]\right) \geq 1- n^{-2},
		$$
		for some constant $C_{\est}>0$ independent of $\alpha, n, d,\beta,\betamax,\betamin$, where $n_0$ depends only on $b$, $l_0$ and the constants of the class.
	\end{Proposition}
	Proposition \ref{prop:smooth-est} reveals that our $\beta_{\est} $ is a good ``undersmooth'' estimator for $\beta$, which ensures that with high probability, we have $\cH(\beta,L)\subseteq \cH(\beta_{\est},L')$ for a constant $L'$ depending on $L,\betamin,\betamax$. This is important for the subsequent plug-in approach. Based on the smoothness estimation, we present our smoothness adaptive algorithm in Algorithm \ref{alg:adaptive-bandit}, with the corresponding regret control in Theorem \ref{thm:adaptive-regret-ub}.
	\begin{algorithm}[h!]
		\caption{Smoothness Adaptive Single-Index Bandits}
		\label{alg:adaptive-bandit}
		\begin{algorithmic}
			\REQUIRE $\betamin$, $\betamax$, $N_0$
			\FOR{$t=1,...,2KN_0$}
			\STATE  Run the smoothness estimator and get $\beta_{\est}$.
			\ENDFOR

			\FOR{$t=2KN_0+1,...,n$}
			\STATE  Run Algorithm \ref{alg:sim-bandit} with smoothness level $\beta_{\est}$ with sample index $t'\gets t-2KN_0$.
			\ENDFOR

			\STATE \textbf{return} $\{\pi_t\}_{t=1}^n$
		\end{algorithmic}
	\end{algorithm}

	\begin{Theorem}\label{thm:adaptive-regret-ub}
		Under the conditions of Proposition \ref{prop:smooth-est}, suppose when run the smoothness estimator for the first $2KN_0$ rounds, and run Algorithm \ref{alg:sim-bandit} for the rest with smoothness level $\beta_{\est}$, then we have

		\begin{equation*}
			\sup_{\cP,\cH^{ss}\left(\beta,L, b, l_0\right)} \Reg(\pi)    \leq C (\log n)^{^{ \iota_2(\alpha,\beta,\betamax,\betamin) }} \cdot \left[ n\left(\frac{ 1 }{n}\right)^{\frac{ \beta (1+\alpha) }{2\beta+1 } }  \vee 1\right],
		\end{equation*}
		provided that $ \alpha(\beta\wedge 1)\le 1$, $\betamin\betamax\ge1$, $\betamax>\frac{2\betamin-1}{2(1-\betamin)}$, $d \lesssim   n^{ \frac{ 2(1-\betamin)\betamax+1-2\betamin}{(2\betamax+1)^2 }  }\wedge n^{\frac{2\betamin-1}{2\betamin+1}}$, and Assumption~\ref{asp:decision-region:2} holds with $c_{\sfm}n^{-\frac{1}{2\betamin+1}}$ in place of $c_{\sfm}n^{-\frac{1}{2\beta+1}}$, for some constant $C$ that depends only on the constants defining $\cP$, on $\beta,L,\alpha,\betamin,\betamax$ and on $K/\lp_K$, so that it is free of $n$ and depends on $d$ only through the constants of $\cP$. Here $\iota_2$ depends only on $\alpha,\beta,\betamax,\betamin,C_{\est}$, and $\iota_2= \frac{1}{2\betamin}+\frac{\beta(1+\alpha)}{2\beta+1} + C_{\est} (1+\alpha) \frac{(2\betamax+1)^2}{\betamin(2\betamin+1)^2} $ when $\beta+1-\alpha\beta > 0$.
	\end{Theorem}

	\begin{Remark}[Smoothness range and exploration cost]\label{rmk:adapt-cost} The $2KN_0$ exploration rounds are at most of the order of the optimal regret for every $\beta$ in the range when $N_0\lesssim n^{\frac{1}{2\betamax+1}}$, which is the first condition on $d$ in Theorem~\ref{thm:adaptive-regret-ub}. Its exponent is positive exactly when $\betamax>\frac{2\betamin-1}{2(1-\betamin)}$, so that $\betamin=0.9$ requires $\betamax>4$; the range $(0.9,1.9)$ of our simulations is not covered.
	\end{Remark}

	\section{Single-Index Bandits in Increasing Dimensions}
	\label{sec:increasing-dim}
	A large body of the literature on single-index models focuses on the case where the dimension $d$ is fixed with respect to $n$. However, in many practical problems, the dimension may increase with $n$ (while still satisfying $d \ll n$ for learnability). In this section, we study the fundamental difficulty of single-index bandits in increasing-dimensional settings, where $d$ can grow with $n$ but remains much smaller than $n$.

Since $d$ now varies, we require the design to be normalized so that the constants of the distribution class do not vary with it.
\begin{Assumption}[Design in increasing dimensions]\label{asp:normalization} For each $d$, let $\cP^{(d)}_0$ be the class of Definition~\ref{def:dist-bandit}, with the boundedness condition~\eqref{eq:design-bounded} replaced by uniformly sub-Gaussian projections, $\sup_{u\in\bbB_2(1)}\norm{u^\top X}_{\psi_2}\le C_{\psi}$. The constants defining $\cP^{(d)}_0$, namely those in Assumptions~\ref{asp:model}, \ref{asp:sufficient-cond}, \ref{asp:nondeg}, \ref{asp:decision-region} and~\ref{asp:margin}, and in particular $C_{\psi}$ and $c_{\sfX}$, are free of $d$. The class $\cP^{(d)}$ consists of the bandit distributions $P$ that share with some $Q\in\cP^{(d)}_0$ a conditional law of the rewards given $X$ satisfying Assumption~\ref{asp:model} at every context, and have $\norm{P_X-Q_X}_{\mathrm{TV}}:=\sup_{A}\abs{P_X(A)-Q_X(A)}\le n^{-2}$.
\end{Assumption}
Since the constants in Theorems~\ref{thm:est}, \ref{thm:minimax-ub} and~\ref{thm:adaptive-regret-ub} are functions of the constants in those assumptions, Assumption~\ref{asp:normalization} makes them free of $d$ as well, so that the rates below are sharp in $d$. The tolerance $n^{-2}$ changes the regret of any policy by at most $1$, so the minimax regrets over $\cP^{(d)}$ and over $\cP^{(d)}_0$ differ by at most $1$ and the rates below hold over either class; it serves only to let the hard instances of Theorem~\ref{thm:minimax-lb-phase} use designs on which the density conditions fail on a set of probability at most $n^{-2}$. No bound on the index range is needed, since Assumption~\ref{asp:sufficient-cond:2} gives $\Leb(I_0)\le1/\lmu$. It is satisfied by coordinate-bounded designs with $\Cov(\widetilde{X})\succeq c_{\sfX}I_{d-1}$, the analogue of the covariate diversity condition $\lambda_{\min}(\Cov(X))\ge c$ in the linear bandit literature with which we compare \citep{bastani2020online,bastani2021mostly,mao2025optimal}. Section~D.3 of the Supplementary Material~\citep{suppSIB} discusses these designs, proves the bound on the effect of the tolerance, and explains how the index error is transferred to the reward under this normalization.
	\subsection{Phase transition in increasing dimensions}\label{sec:phase-trans-1}

	Interestingly, in the case when $d$ is increasing, the minimax rate exhibits a phase transition between the nonparametric and linear regimes. Specifically, when $d$ is small, the optimal rate is governed by the nonparametric phase, where the smoothness parameter $\beta$ determines the regret rate. In contrast, when $d$ becomes large, the optimal rate enters the parametric phase, as learning the index becomes more challenging; in this case, the optimal regret corresponds to that of linear bandits \citep{bastani2021mostly, mao2025optimal}, which is independent of $\beta$.

	We show that in the increasing dimensions case, our Algorithm \ref{alg:sim-bandit} shares a two-phase regret upper bound. For simplicity, we assume $\beta\ge 1$ and is known.

	\begin{Theorem}[Regret upper bound]\label{thm:minimax-ub-phase} Suppose $\beta \ge 1$ and $d\ll n$, and let Algorithm~\ref{alg:sim-bandit}, denoted by $\widehat{\pi}$, use $(d+\log^2 n)\log n$ in place of $d+\log^2 n$ in $n_m$. Then it achieves the regret
		$$
		\sup_{\cP^{(d)},\cH(\beta,L)}	\Reg(\widehat{\pi})   \leq C (\log n)^{ \iota_3(\alpha,\beta) } \cdot \left[ n\left(\left(\frac{ 1 }{n}\right)^{\frac{ \beta (1+\alpha) }{2\beta+1 } } +   \left(\frac{d}{n}\right)^{\frac{ 1+\alpha}{2}\wedge 1}\right)\vee 1\right],
		$$
		where, under Assumption~\ref{asp:normalization}, $C>0$ depends only on the constants defining $\cP^{(d)}$ and on $\beta,L,\alpha$, is free of $n$ and $d$, and is polynomial in $K/\lp_K$ as in Theorem~\ref{thm:minimax-ub}, and $\iota_3=\iota_1\vee 4$ with $\iota_1$ as in Theorem~\ref{thm:minimax-ub}. Specifically, when $\alpha\le1$,
		\begin{equation}\label{eq:regret-ub-phase}
			\sup_{\cP^{(d)},\cH(\beta,L)}	\Reg(\widehat{\pi})   \lesssim  (\log n)^{ \iota_3(\alpha,\beta) }\cdot\begin{cases} n^{1-\frac{(1+\alpha)\beta}{2\beta+1}} & \text { if } d\lesssim n^{\frac{1}{2\beta+1}},
				\\
				n \left(\frac{d}{n}\right)^{\frac{1+\alpha}{2}} & \text { if }d\gtrsim n^{\frac{1}{2\beta+1}}. \end{cases}
		\end{equation}

	\end{Theorem}
	Notably, our algorithm achieves this two-phase regret bound without prior knowledge of the phase. The first rate is that of the nonparametric phase (Section~\ref{sec:opt-rates}) and the second that of linear bandits; when $\alpha > 1$ the latter becomes $\widetilde{O}(d)$, the cost of learning in $\R^d$.
	\begin{Theorem}[Regret lower bound]\label{thm:minimax-lb-phase} Under Assumption~\ref{asp:normalization}, suppose that $\beta\ge1$, and that $\alpha\beta\le1$ or $\alpha=1$. Then for any admissible policy $\pi$,
		\begin{equation}\label{eq:minimax-lb-phase}
			\inf _\pi \sup _{\cP^{(d)},\cH(\beta,L)} 	\Reg(\pi)   \ge c\cdot \begin{cases} n^{1-\frac{(1+\alpha)\beta}{2\beta+1}} \vee 1 & \text { if } d\lesssim n^{\frac{1}{2\beta+1}},
				\\
				n \left(\frac{d}{n}\right)^{\frac{1+\alpha}{2}} & \text { if }n^{\frac{1}{2\beta+1}}\lesssim d\lesssim \frac{n}{\log n}, \end{cases}
		\end{equation}
		where $c>0$ is a constant independent of $n$ and $d$.
	\end{Theorem}
	Theorem~\ref{thm:minimax-lb-phase} shows that, under its conditions, the two-phase rate of Theorem~\ref{thm:minimax-ub-phase} is minimax optimal up to logarithmic factors. For $\beta\ge1$ its condition coincides with that of Theorem~\ref{thm:minimax-lb}. It also implies $\alpha\le1$, which is the regime in which the parametric term in~\eqref{eq:regret-ub-phase} takes the form $n(d/n)^{\frac{1+\alpha}{2}}$; the case $\alpha>1$ remains open even for linear bandits \citep{goldenshluger2009woodroofe, li2021regret}.
	In this range the minimax regret is therefore of order $n^{1-\frac{(1+\alpha)\beta}{2\beta+1}}+n(d/n)^{\frac{1+\alpha}{2}}$, and single-index bandits interpolate between one-dimensional nonparametric bandits and linear bandits, with the transition at $d\asymp n^{\frac{1}{2\beta+1}}$. Below this threshold the difficulty lies in learning the links, and above it in learning the indices.

	We briefly explain the hard instance construction for the parametric phase, which is more challenging than the nonparametric one. When $\alpha=1$ or $\beta=1$, we construct two arms sharing one index, whose link difference changes sign about $(d/n)^{\frac{\alpha-1}{2}}$ times, so that an index perturbation of order $\sqrt{d/n}$ reverses the optimal arm near each sign change. However, when $\beta>1$, a link difference of smoothness $\beta$ cannot change sign this often with a non-vanishing slope, and we instead give the two arms different indices. The perturbation then shifts the level of an oscillating reward difference rather than the positions of its zeros, and moves the decision boundary across a whole half-period, a distance of order $(d/n)^{\frac{1}{2\beta}}$, which is much larger than the perturbation itself.

\section{Numerical Experiments}
\label{sec:simulations}

We report the growth of the regret of Algorithm~\ref{alg:sim-bandit} in the dimension and in the horizon, a comparison with nonparametric, parametric and neural alternatives, and an analysis of four real classification data sets. Section~F of the Supplementary Material~\citep{suppSIB} contains the full settings and the remaining experiments: online index estimation, the adaptive Algorithm~\ref{alg:adaptive-bandit} and the smoothness estimator, the geometry of the margin regions, heavy-tailed rewards (which lie outside Assumption~\ref{asp:model}), the neural comparison in full, the dependence on the number of arms, and timings.

\subsection{Synthetic data simulations}\label{sec:sim-synthetic}

\paragraph*{Set up} We take $K=3$ arms and rewards $Y_t^{(k)}=f_k(v_k^\top X_t)+\xi_t$ with $\xi_t\sim\cN(0,0.1^2)$ independent of $X_t$, and index vectors with $v_{k,1}=1$ and $\norm{\widetilde{v}_k}_2=1$, drawn at random subject to linear independence. Two families of links are used,
\begin{equation}\label{eq:links-scaled}
	f_1(z) = 0.8\,\sgn(z)\abs{\tfrac{z}{2}}^{\beta},\ \
	f_2(z) = 0.5\,\sgn(z)\abs{\tfrac{z}{2}}^{\beta}+ 0.1\,z,\ \
	f_3(z) = 1.5\,\sgn(z)\abs{\tfrac{z}{2}}^{\beta},
\end{equation}
the family of the first version of this paper, whose links differ only by a scale factor, and
\begin{equation}\label{eq:links-curv}
	f_1(z) = 0.8\tanh(2z),\ \
	f_2(z) = 0.15\,z+0.20\tanh(z),\ \
	f_3(z) = \sgn(z)\abs{\tfrac{z}{3}}^{5/2}+0.10\,z ,
\end{equation}
whose links are saturating, nearly linear and strongly convex. The policies of LinUCB and UCB-GLM are the argmax of $K$ affine functions of the context whatever the data, and the best such policy gives up $0.0\%$ of the mean reward gap per round on~\eqref{eq:links-scaled} but $5.6\%$ on~\eqref{eq:links-curv}. At $d=4$ the design is the ball-truncated Gaussian $X_t\sim\cN(0,I_d)\mid \norm{X_t}_2\le 1$ of the first version, and in growing dimension a coordinate-bounded design, which Assumption~\ref{asp:normalization} admits; the measured margin exponent is $\widehat{\alpha}=1.0$ throughout, the generic value for links that cross transversally, and it is the value used when a predicted curve is drawn. The links take values outside $[0,1]$, which Definition~\ref{def:smooth-class} assumes; a common increasing affine map of the rewards brings them into $[0,1]$, leaves the index estimates unchanged, and multiplies the reward differences, the regret and the appropriate $c_{\epsilon}$ by the same factor. The theory asks only that $C_{\cT}$ be large and $c_{\epsilon}$ small relative to constants of the class, and we take $C_{\cT}=1$ and $c_{\epsilon}=0.5$ from a pilot grid. Two further choices are practical devices. The bandwidth constant $C_{\sfH}$, which the theorems treat as fixed, is selected on held-out data within each epoch. An arm is eliminated only if it trails the leader by more than $\epsilon_m$ after the leader's estimate is lowered by half the discrepancy between the two cross-fitting folds of Algorithm~\ref{alg:MRC-reg}; this makes elimination more conservative, and Theorem~\ref{thm:minimax-ub} continues to hold for it with a different constant. Section~F.1 of the Supplementary Material~\citep{suppSIB} gives the grid, this rule, the exploration constant of each competitor and the instance diagnostics. Every panel reports the mean over independent replications and the central $95\%$ of them, with the count in the title.

\begin{Remark}\label{rmk:computation}
	Algorithm~\ref{alg:sim-bandit} requires maximizing the non-convex criterion~\eqref{eq:mrc}, but only once per epoch, and $M\lesssim\log n$, so the total cost remains manageable for large $n$. We optimize it by differential evolution below; the smoothed criterion of Section~\ref{sec:discuss} is $40$ to $100$ times faster at no loss of accuracy.
\end{Remark}

\begin{figure}[t!]
	\centering
	\includegraphics[width=\textwidth]{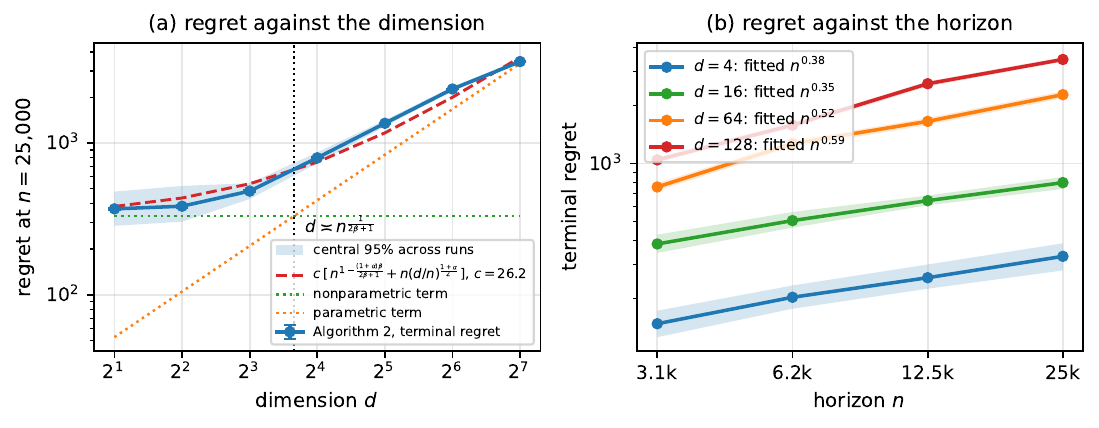}
	\caption{Regret of Algorithm~\ref{alg:sim-bandit} on the links~\eqref{eq:links-scaled}, $\beta=1.5$, coordinate-bounded design, $K=3$. (a) terminal regret against the dimension at $n=25{,}000$ over $200$ replications, with the fitted two-phase curve, its two terms, and the predicted transition $d\asymp n^{1/(2\beta+1)}=12.6$. (b) terminal regret against the horizon at four dimensions over $100$ replications. Shaded regions give the central $95\%$ of the replications, error bars the $95\%$ interval for the mean.}
	\label{fig:rates}
\end{figure}

\paragraph*{Rates in the dimension and the horizon} Figure~\ref{fig:rates} reports the two dependencies of Theorem~\ref{thm:minimax-ub-phase} on the links~\eqref{eq:links-scaled} with $\beta=1.5$. In panel (a) the regret at $n=25{,}000$ is flat over the first three dimensions, $369$, $384$ and $483$ at $d=2,4,8$, and then grows, $801$, $1353$, $2278$ and $3456$ at $d=16,\dots,128$: the ambient dimension does not enter the nonparametric exponent, and enters only through an additive parametric term, negligible until $d\asymp n^{1/(2\beta+1)}=12.6$ and dominant after it. With $\alpha=1$ the bound is $c\,(n^{1/(2\beta+1)}+d)$, and a single constant $c$ fitted on the logarithmic scale reproduces the measured curve over a $64$-fold range of $d$, the implied $c$ varying only between $23$ and $30$ across the seven dimensions. Panel (b) varies the horizon at four dimensions. The fitted exponents in $n$ are $0.38$ and $0.35$ at $d=4$ and $16$, against the predicted $0.25$, and $0.52$ and $0.59$ at $d=64$ and $128$, where the parametric term takes the form $n(d/n)^{(1+\alpha)/2}$ only once the index error $\sqrt{d/n}$, here between $0.20$ and $0.07$, is of second order. Measured exponents exceed predicted ones everywhere, which we attribute to the realized accuracy of $\widehat{g}_{k,m}$, $2.5$ to $5.6$ times coarser than the $\epsilon_m/16$ the good event requires. We therefore read the panels as confirming the shape of the bound and the location of the transition.

\begin{figure}[t!]
	\centering
	\includegraphics[width=\textwidth]{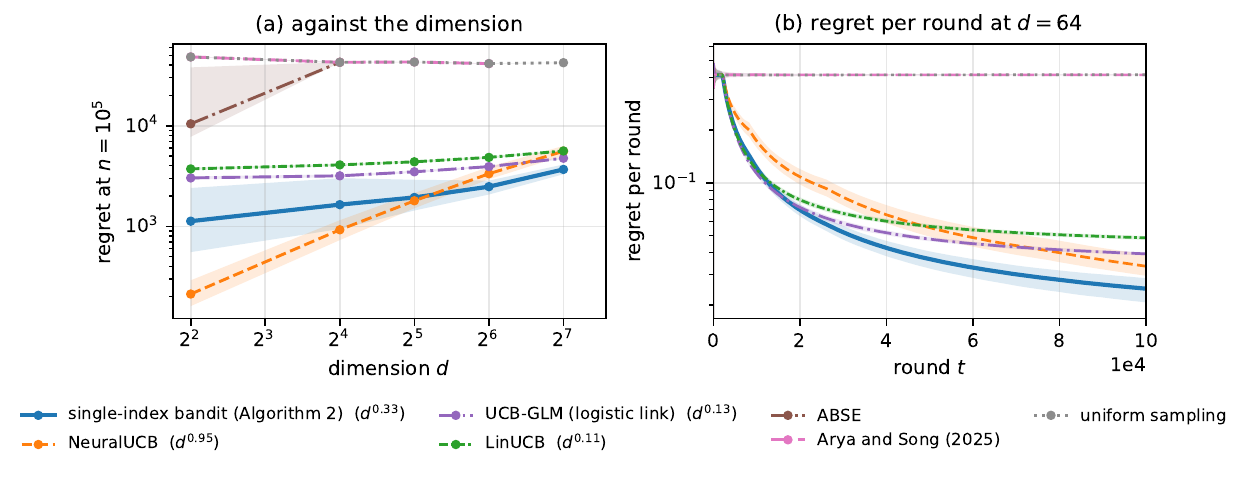}
	\caption{Regret on the links~\eqref{eq:links-curv}, coordinate-bounded design, $K=3$, $n=10^5$, $100$ replications, with every method held to the $M=O(\log n)$ policy updates of Algorithm~\ref{alg:sim-bandit}. (a) terminal regret against the dimension, with the fitted power of $d$ in the legend. (b) regret per round at $d=64$ on a logarithmic scale, so that all six procedures that can be run at that dimension fit in one panel. Shaded regions give the central $95\%$ of the replications. ABSE refines a dyadic partition of $\cX$ and is shown only where that partition can still be built, $d\le 16$; the procedure of \cite{arya2025semi} is shown up to $d=64$ and lies on the uniform-sampling curve at every dimension; SmoothBandit is omitted, a single replication costing more than five minutes already at $d=4$ and this horizon.}
	\label{fig:comparison}
\end{figure}

\paragraph*{Comparison with alternative procedures} Figure~\ref{fig:comparison} compares Algorithm~\ref{alg:sim-bandit} on the links~\eqref{eq:links-curv} at $n=10^5$ over $100$ replications with the adaptively binned successive elimination algorithm of \cite{perchet2013multi} (ABSE), the batched single-index procedure of Arya and Song \citep{arya2025semi}, LinUCB \citep{li2010contextual,chu2011contextual}, UCB-GLM \citep{li2017provably}, NeuralUCB \citep{zhou2020neural} and uniform sampling. Every method is held to the epoch schedule of Algorithm~\ref{alg:sim-bandit}, so that all of them update their policy $M=O(\log n)$ times; this is the protocol the theory of this paper describes. ABSE and the smooth-bin successive elimination algorithm of \cite{hu2022smooth,gur2022smoothness} both partition $\cX$ and cannot be run much past $d=16$; at $d=4$ our advantage over them is a factor between $5.7$ and $9.0$ (Section~F.3 of the Supplementary Material~\citep{suppSIB}).

Panel (a) varies the dimension. The regret of Algorithm~\ref{alg:sim-bandit} grows like $d^{0.33}$, since the dimension enters our bound only through the additive parametric term, while NeuralUCB's grows like $d^{0.95}$; the two curves cross between $d=32$ and $d=64$. LinUCB and UCB-GLM are nearly flat in $d$ but sit well above both, because what limits them is the bias of an affine-argmax policy and not estimation error. ABSE pays the exponential price in the dimension that the single-index structure avoids. The procedure of Arya and Song coincides with uniform sampling, for a reason specific to these instances. Their model $g_k(x)=f(v^\top x)+\theta_k$ ranks the arms the same way at every context, so their elimination test compares the arms' \emph{marginal} means, and we have built the instance so that those are equal. Nothing is ever eliminated. The fault is ours and not theirs: on their own model the same implementation reduces to a single arm and beats uniform sampling by two orders of magnitude (Section~F.4 of the Supplementary Material~\citep{suppSIB}). Panel (b) gives the regret per round at $d=64$, where LinUCB and UCB-GLM level off at their bias floors while ours is still falling at the end of the horizon.

We also add two qualifications, both developed in Section~F.5 of the Supplementary Material~\citep{suppSIB}. At $d\le 32$, and at every $d$ if the baselines are left in their own fully sequential form, NeuralUCB is ahead. The difference is not one of estimation, since Algorithm~\ref{alg:MRC-reg} is the more accurate of the two in many settings. It is one of exploration. NeuralUCB is effectively greedy, and it borrows the exploration it needs from the covariate diversity that Assumption~\ref{asp:nondeg} provides; Algorithm~\ref{alg:sim-bandit} has to pay for that exploration, because its guarantee must hold when the diversity is absent.

\subsection{Real data analysis}\label{sec:sim-real}

We next apply the method to four real data sets: the UCI mushroom data ($8124$ instances, $22$ categorical attributes one-hot encoded, $2$ classes) and the Statlog shuttle data ($9$ attributes, the four most frequent of its $7$ classes, $57{,}927$ instances) from OpenML, and Wisconsin diagnostic breast cancer ($569$, $2$) and wine ($178$, $3$) as distributed with \texttt{scikit-learn} \citep{pedregosa2011scikit}. Features are standardized and reduced to $d=8$ by principal components; the arms are the classes, a context is drawn at random from the pool each round, and the reward is $1$ if the selected arm is the true class, so the regret counts misclassifications. We use $n=10^5$ and $20$ replications. Here the single-index structure is an assumption about the data rather than something we impose, so we first examine it diagnostically.

\begin{figure}[t!]
	\centering
	\includegraphics[width=\textwidth]{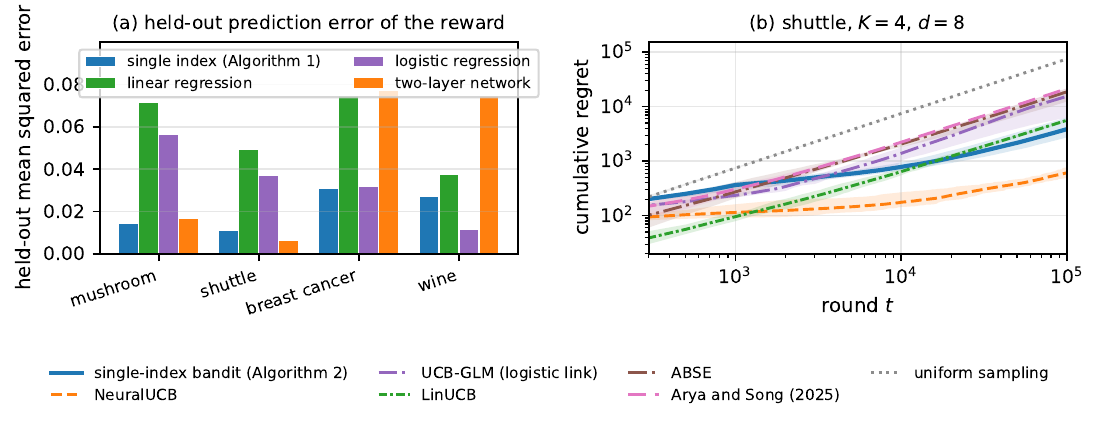}
	\caption{Real data, $d=8$ after principal components. (a) held-out mean squared error of the fitted reward function, averaged over the arms of each data set, for Algorithm~\ref{alg:MRC-reg} and three alternatives. (b) cumulative regret, equal to the number of misclassifications, of seven procedures on the Statlog shuttle data at $n=10^5$ over $20$ replications, on a logarithmic scale in both axes, with the central $95\%$ of the replications shaded. SmoothBandit is omitted, its bin schedule being built for $K=3$; the procedure of \cite{arya2025semi}, whose ranking of the arms cannot depend on the context, settles on the most frequent class.}
	\label{fig:realdata}
\end{figure}

Panel (a) of Figure~\ref{fig:realdata} fits each arm's reward function on a random $70\%$ of the pool and reports the mean squared error on the rest; all tuning, including the bandwidth, uses the $70\%$ only, so that the remaining $30\%$ serves for evaluation alone. The single-index fit attains $0.014$, $0.011$, $0.031$ and $0.027$ on mushroom, shuttle, breast cancer and wine, against $0.072$, $0.049$, $0.075$ and $0.037$ for linear regression and $0.017$, $0.006$, $0.077$ and $0.075$ for a network with one hidden layer of $32$ units. It beats the linear model everywhere, and the network everywhere but shuttle. Logistic regression, which is itself a single-index model with a known link, is comparable on breast cancer, better on wine, and worse on mushroom and shuttle. The held-out rank correlation between the reward and the fitted index is also positive for all $11$ arms. These diagnostics are consistent with a single-index model with rank-monotone responses, but they do not establish it: a positive aggregate correlation does not imply monotonicity at the level of individual contexts, and the smoothness and margin conditions of the regret theorems are not checked. We therefore present these data as a practical benchmark rather than as a validation of the theory.

Panel (b) turns to the bandit, on shuttle, where the regret is the number of misclassifications. At $n=10^5$ Algorithm~\ref{alg:sim-bandit} ends at $3865$ $[2721,6351]$, against $5628$ for LinUCB, $15{,}449$ for UCB-GLM, $18{,}765$ for ABSE and $21{,}365$ for the procedure of Arya and Song \citep{arya2025semi}; only NeuralUCB, at $610$, does better. Our procedure therefore beats every fully nonparametric competitor, the nearby single-index procedure, and both parametric methods with a prescribed link. Mushroom gives the same ordering, and there the gap to LinUCB grows with the horizon, from a factor $1.4$ at $n=2\cdot10^4$ to $1.8$ at $n=10^5$ (Section~F.8 of the Supplementary Material~\citep{suppSIB}).

Panel (a) shows our reward estimates to be more accurate than LinUCB's on all four sets and than NeuralUCB's on three, so the remaining gap to NeuralUCB may be one of exploration rather than estimation: it runs almost greedily here, the exploration it needs being supplied by the covariate distribution itself, the regime in which greedy policies are known to be preferred \citep{bastani2021mostly}. That is a property of the instance and not a guarantee uniformly for all data.

\paragraph*{Further experiments} We examine three further questions. We track the index estimate across the epochs of Algorithm~\ref{alg:sim-bandit} and measure its offline error: it improves monotonically over the epochs, and it scales as $\sqrt{d/n}$ up to a constant, the rate of Theorem~\ref{thm:est}. We compare the adaptive Algorithm~\ref{alg:adaptive-bandit} with Algorithm~\ref{alg:sim-bandit} given the true smoothness. The two perform comparably, and at low noise the estimated smoothness is ordered in the truth and falls below it, which is the undersmoothing Proposition~\ref{prop:smooth-est} calls for. The probe rounds are expensive, though: at the sample sizes we can run they take up a large share of the horizon and account for most of the regret. Finally, we vary the margin exponent by construction. The fitted regret exponent falls as $\alpha$ grows, as $1-(1+\alpha)\beta/(2\beta+1)$ requires. We also build an instance whose optimal regions break into many thin slivers, which violates the regularity of Assumption~\ref{asp:decision-region:2} while leaving $\alpha$ fixed; it costs almost nothing in regret, so the requirement looks like something our proof needs rather than something the algorithm does. Sections~F.2, F.6 and~F.7 of the Supplementary Material~\citep{suppSIB} report these in full.

\section{Discussion}
\label{sec:discuss}

\textit{Open problems.} We briefly mention four directions. (i) \textit{Non-monotone and multiple-index rewards.} Relaxing rank monotonicity would likely call for a different index estimator under bandit sampling (Remark~\ref{rmk:monotone}); the multiple-index case $g_k(x)=f_k(V_k^\top x)$ with $V_k\in\R^{d\times q}$ would interpolate between our rates and the $q$-dimensional nonparametric rates. (ii) \textit{Sparse indices.} When $\norm{v_k}_0\le s\ll d$, one expects the parametric phase to improve to $n(s\log d/n)^{\frac{1+\alpha}{2}}$, with the transition at $s\log d\asymp n^{\frac{1}{2\beta+1}}$; the natural estimator is the $\ell_1$-constrained rank correlation, for which the curvature identity should yield a restricted-eigenvalue version of~\eqref{eq:curvature-main}. (iii) \textit{Computation.} The rank correlation criterion is non-convex and costs $O(n_m^2d)$ per epoch. Since our bounds are driven by Hoeffding decompositions of the associated $U$-process, an incomplete $U$-statistic built from $O(n_m)$ random pairs, or a smoothed criterion in the spirit of \cite{horowitz1992smoothed}, should preserve the rates at cost $O(n_md)$ while permitting gradient-based optimization. (iv) \textit{The remaining range of the lower bound.} Theorem~\ref{thm:minimax-lb} matches Theorem~\ref{thm:minimax-ub} for $\alpha\beta\le1$ and for $\alpha=1\le\beta$, but not yet for $\beta>1$ and $1/\beta<\alpha<1$, although this range lies within $\alpha(\beta\wedge1)\le1$, where the upper bound takes its simplest form.

\begin{funding}
The research  was supported in part by NSF Grant DMS-2413106 and NIH grants R01-GM129781 and R01-GM123056.
\end{funding}

\begin{supplement}
\stitle{Supplementary Materials to “Nonparametric Bandits with Single-Index Rewards:
	Optimality and Adaptivity”}
\sdescription{The supplementary file contains the proofs of all results in the main text (the separation and curvature lemma for the rank correlation, the non-asymptotic single-index regression bound, the upper and lower bounds for the regret, and the adaptivity results), together with the threshold for a static oracle policy, a verification of the assumptions for two canonical families of instances, an expanded review of the related literature, and the remaining numerical results: the tuning grids, the full comparisons with competing procedures, heavy-tailed rewards, the margin geometry, the real data analysis, and timings.}
\end{supplement}

\bibliographystyle{imsart-number}
\bibliography{reference}

\end{document}